\documentclass[11pt, leqno, twoside]{article}
\usepackage[left=3cm,top=3cm,right=3cm,bottom=4.5cm]{geometry}
\usepackage{amssymb, amsmath} 
\usepackage{epic} 
\usepackage{tikz}
\usetikzlibrary{trees}
\usepackage{enumerate}
\usepackage{color}
\usepackage{mathtools}
\usepackage{amsthm}
\newtheorem{thm}{Theorem}%[section]
\newtheorem{cor}[thm]{Corollary}
\newtheorem{lem}[thm]{Lemma}
\newtheorem{prop}[thm]{Proposition}
\newtheorem{defn}[thm]{Definition}

\DeclarePairedDelimiter\ceil{\lceil}{\rceil}
\DeclarePairedDelimiter\floor{\lfloor}{\rfloor}

\newcommand{\N}{\mathbb{N}}
\newcommand{\R}{\mathbb{R}}

\newcommand{\ve}{\varepsilon}

\newcommand{\Om}{\Omega}

\DeclareMathOperator{\diag}{diag}

\DeclareMathOperator{\rank}{rank}
\DeclareMathOperator{\supp}{supp}

\DeclareMathOperator{\id}{id}

\tikzstyle{level 1}=[level distance=4cm, sibling distance=2.5cm]
\tikzstyle{level 2}=[level distance=3cm, sibling distance=2cm]

\tikzstyle{bag} = [text width=10em, text centered]
\tikzstyle{end} = [circle, minimum width=3pt,fill, inner sep=0pt]

\begin{document}
%\pagestyle{myheadings} 
%\markboth{\texttt{\jobname{.tex}}}{\texttt{\jobname{.tex}}}

\title{Bi-Sobolev homeomorphisms $f$ with $Df$ and $Df^{-1}$ of low rank using laminates}

\author{Marcos Oliva\\
\footnotesize Department of Mathematics, Faculty of Sciences, Universidad Aut\'onoma de Madrid, E-28049 Madrid, Spain\\
\footnotesize ICMAT CSIC-UAM-UCM-UC3M, E-28049 Madrid, Spain}

\date{\today}
\maketitle
\begin{abstract}
Let $\Omega\subset \R^{n}$ be a bounded open set. Given $1\leq m_1,m_2\leq n-2$, we construct a homeomorphism $f :\Omega\to \Omega$ that is H\"older continuous, $f$ is the identity on $\partial \Omega$, the derivative $D f$ has rank $m_1$ a.e.\ in $\Omega$, the derivative $D f^{-1}$ of the inverse has rank $m_2$ a.e.\ in $\Omega$, $Df\in W^{1,p}$ and $Df^{-1}\in W^{1,q}$ for $p<\min\{m_1+1,n-m_2\}$, $q<\min\{m_2+1,n-m_1\}$. The proof is based on convex integration and laminates. We also show that the integrability of the function and the inverse is sharp.
\end{abstract}
\section{Introduction}\label{sect: introduction}
In this paper we prove that there exists a bi-Sobolev homeomorphism $f:\Omega\to\R^{n}$ whose derivative $Df$ has rank equal to $m_1$ at a.e.\ point and $\rank (Df^{-1}(y))=m_2$ for a.e. $y\in f(\Omega)$, for given $m_1,m_2 \leq n-2$. This map cannot satisfy the Luzin condition (which states that the image by $f$ of a set of measure zero also has measure zero).

There are several examples in the literature of Sobolev homeomorphisms that do not satisfy the Luzin condition, see  \cite{Cerny11}, \cite{Hencl11}. See, in addition, \cite{KaKoMa01} for some sufficient conditions for its validity.

The pioneering construction by Hencl \cite{Hencl11} of Sobolev homeomorphisms with zero Jacobian has been further developed in \cite{Cerny15,DoHeSc14} to construct bi-Sobolev homeomorphisms $f$ with zero minors of $Df$ and of $Df^{-1}$ from some order. See \cite{FaMoOl16} for further state of the art.

As far as we know, there are two methods to obtain this class of pathological maps:
\begin{itemize}
\item The method showed in \cite{Hencl11} and its developments \cite{DoHeSc14,Cerny15}. All those constructions were based on a careful explicit construction and a limit process to obtain a Cantor set where the Jacobian is supported. 
\item The method that appears in \cite{FaMoOl16}, based on staircase laminates. The use of laminates in relation with integrability issues was initiated in \cite{Faraco03}, where the so-called staircase laminates were introduced. 
This technique has turned out to be extremely efficient in a number of unrelated results, such as \cite{Pedregal93}, \cite{Muller96}, \cite{CoFaMa05}, \cite{CoFaMaMu05}, \cite{AsFaSz08} and \cite{BoSzVo13}; see also \cite{KiKr11} and \cite{KiKr16}, where rank-one convexity is used to obtain badly integrable maps. In \cite{AsFaSz08} a version of the in-approximation of M\"uller and \v{S}ver\'ak \cite{MuSv03} valid for unbounded set was developed. In \cite{FaMoOl16}, these ideas were extended to construct extremal homeomorphisms: this was an extra feature that did not appear before.
\end{itemize}
However, in none of those papers there was a consideration on the inverse map.
In particular, since the laminates behind the proof of \cite{FaMoOl16} were supported in a set of non-invertible matrices, here we need to proceed differently. Let $E_m$ be the set of symmetric matrices of rank $m$. Then we construct a suitable sequence of open sets $E^j_m$ of symmetric positive definite matrices converging to $E_m$. Then, when we denote the inverse of a laminate $\nu=\sum_{i=1}^{N}\lambda_i\delta_{A_i}$ supported in the set of invertible matrices as 
\[\nu^{-1}=\frac{1}{\det(\sum_{i=1}^{N}\lambda_i A_i)}\sum_{i=1}^{N}\lambda_i\det (A_i)\delta_{A_i^{-1}}\]
(see Definition \ref{de:inverse laminate}), we construct a sequence of staircase laminates $\nu_j$ supported in $E^j_{m_1}$ such that the inverse laminates $\nu_j^{-1}$ are supported in $E_{m_2}^{j}$. We run convex integration to find sequences $f_j$ such that the gradient of $f_j$ is arbitrarily close to $\nu_j$, and the gradient of $f_j^{-1}$ is arbitrarily close to $\nu_j^{-1}$.
Then, following the ideas of \cite{MuSv03} and \cite{AsFaSz08} we are able to show that both $f_j$ and $f_j^{-1}$ converge to maps $f$ and $f^{-1}$ whose gradients are supported in $E_{m_1}$ and $E_{m_2}$, respectively, and have the highest possible integrability for a bi-Sobolev homeomorphism with rank of the gradient equal to $m_1$ and rank of the gradient of its inverse equal to $m_2$. This method allows us to obtain a pathological bi-Sobolev map, instead of obtaining just a Sobolev map.

Thus, we show that with laminates it is possible to combine the results of \cite{AlAm99,Faraco06,Hencl11}. In fact, we also improve the result
of \v{C}ern\'{y} \cite{Cerny15}, where he constructs a bi-Sobolev homeomorphism $f$ with all the minors of  order $m_1 +1$ of $Df$, and all the minors of order $m_2 +1$ of $Df^{-1}$ equal to zero almost everywhere, belonging to $W^{1,p}$ with $1\leq p <\min\{\frac{n}{n-m_1},n-m_2\}$ and with the inverse $f^{-1}$ in $W^{1,q}$ with $1\leq q<\min\{\frac{n}{n-m_2},n-m_1\}$. 
The bi-Sobolev homeomorphism $f$ that we construct satisfies that $Df$ has rank equal to $m_1$ and $Df^{-1}$ has rank equal to $m_2$, and $f$ is in $W^{1,p}$ and $f^{-1}$ is in $W^{1,q}$ for $1\leq p<\min\{m_1+1,n-m_2\}$ and $1\leq q<\min\{m_2+1,n-m_1\}$.
We also recover a slightly weaker version of the result in \cite{FaMoOl16}; to be precise, here we obtain that $Df\in L^{p}$ for all $p<\rank(Df)+1$, whereas in \cite{FaMoOl16} it is obtained that $Df$ is in the weak  $L^{\rank(Df)+1}$ space.
 
\begin{thm}\label{theorem existence}
Let $\Omega\subset \R^{n}$ be a bounded open set, $m_1,m_2\in\N$, $1\leq m_1,m_2\leq n-2$, $\ve>0$ and $\alpha\in (0,1)$. Then there exists a convex function $u:\Omega\to\R$ whose gradient $f=\nabla u:\Omega\to\Omega$ is a bi-Sobolev homeomorphism and satisfies:
\begin{enumerate}[i)]
\item\label{th ex border} $f=\id$ on $\partial \Om$.
\item\label{th ex holder} $\|f-\id\|_{C^{\alpha}(\overline{\Omega})}<\ve$ and $\|f^{-1}-\id\|_{C^{\alpha}(\overline{\Omega})}<\ve$.
\item\label{th ex rank} $\rank (D f(x))=m_1$ a.e. $x\in\Omega$ and $\rank (D f^{-1}(y))=m_2$ for a.e.\ $y\in\Om$.
\item\label{th ex int} $D f\in \bigcap_{p<m_1'+1}L^{p}\left(\Omega,\R^{n\times n}\right)$ and $D f^{-1}\in \bigcap_{q<m_2'+1}L^{q}\left(\Omega,\R^{n\times n}\right)$,
\end{enumerate}
where for $i=1,2$ we have
\begin{equation*}
m_i'=\begin{cases}
m_i&\text{if }m_1+m_2\leq n-1,\\
m_i-(m_1+m_2-n+1)&\text{if }m_1+m_2\geq n.
\end{cases}
\end{equation*}
\end{thm}

 Notice that, as explained in \cite{Hencl11}, using the area formula for Sobolev mappings (\cite{Hajlasz93}) we have that this kind of homeomorphisms sends a set of full measure to a null set, and a null set to a set of full measure.
 
 The key to see that Theorem \ref{theorem existence} is sharp is the following result, whose proof is based in the one of \cite[Th. 4]{HeMoPaSb09} and, in fact, generalizes that result.
\begin{thm}\label{theorem integrability of f and rank of the inverse}
Let $f:\Omega\to\R^{n}$ be a continuous bi-Sobolev map such that $f^{-1}\in W^{1,q}\left(f(\Omega),\R^{n}\right)$ and for a measurable set $E\subset\Omega$ we have $Jf=0$ almost everywhere on $E$. Let $m\in \N$ and assume one of the following:
\begin{enumerate}[a)]
\item\label{optionA} $q = m$ and $f^{-1}$ is H\"{o}lder continuous.
\item\label{optionB} $q>m$.
\end{enumerate}
Then $\rank\left(Df\right)\leq n-m-1$ almost everywhere on $E$.
\end{thm}
Observe that this implies that given a continuous bi-Sobolev map $f\in W^{1,p}$ with $p>n-1$, we have $Jf^{-1}(y)\neq 0$ a.e. $y\in f(\Omega)$.

Using the last theorem and Theorem 12 of \cite{FaMoOl16} we obtain the following two results. We think that the next one has interest by itself: it gives an upper bound for the sum of the integrabilities of $Df$ and $Df^{-1}$ for a bi-Sobolev homeomorphism with Jacobian equal to zero almost everywhere.
\begin{thm}\label{theorem sum of the integrabilities}
Let $f:\Omega\to\R^{n}$ be a bi-Sobolev homeomorphism such that $f\in W^{1,p}\left(\Omega,\R^{n}\right)$, $f^{-1}\in W^{1,q}\left(f(\Omega),\R^{n}\right)$, $f|_{\partial\Omega}=\id|_{\partial\Omega}$ and $Jf=0$ almost everywhere in $\Omega$. Then $p+q\leq n+1$. Moreover, if $p+q=n+1$ then $p,q\in\N$ and $f,f^{-1}\notin C^{\alpha}$ for any $\alpha\in (0,1)$.
\end{thm}

We do not know whether there exists a bi-Sobolev homeomorphism that attains the equality $p+q=n+1$.

Finally, the next theorem shows that Theorem \ref{theorem existence} is sharp and gives bounds for the integrability of $Df$ and $Df^{-1}$ depending on the rank of $Df$ and $Df^{-1}$. 
\begin{thm}\label{theorem integrabilities and rank}
Let $m_1,m_2\in\N$ be such that $m_1,m_2\leq n-1$, let $f:\Omega\to\R^{n}$ be a bi-Sobolev homeomorphism such that $f\in W^{1,p}\left(\Omega,\R^{n}\right)$, $f^{-1}\in W^{1,q}\left(f(\Omega),\R^{n}\right)$, $f|_{\partial\Omega}=\id|_{\partial\Omega}$, $\rank(Df)\leq m_1$, $\rank(Df^{-1})\leq m_2$ almost everywhere in $\Omega$ and suppose that there exist measurable $A,B\subset\Omega$ with $|A|,|B|>0$, $\rank(Df)=m_1$ on $A$ and $\rank(Df^{-1})=m_2$ on $B$. Then $p\leq \min\{m_1+1, n-m_2\}$ and $q\leq \min\{m_2+1, n-m_1\}$. 
\end{thm}
Observe that this theorem implies that $p+q\leq n+1$.

The structure of the paper is as follows. Section \ref{sect: notation} presents the general notation of the paper. In Section \ref{sect: sharp} we prove Theorems \ref{theorem integrability of f and rank of the inverse}, \ref{theorem sum of the integrabilities} and \ref{theorem integrabilities and rank}. Section \ref{sect: sketch proof} shows the construction of the laminates in the simplest case of dimension three. This section gives an idea of how we prove Theorem \ref{theorem existence} and helps to understand the actual proof. Section \ref{sect: proof of existence}, which is the bulk of the paper, constructs a sequence of laminates that converges to the probability measure sketched in Section \ref{sect: sharp}, but in general dimensions, as well as a sequence of functions whose gradients approximate the laminates.
\section{General notation}\label{sect: notation}

We explain the general notation used throughout the paper, most of which is standard.

In the whole paper, $\Om$ is an open, non-empty bounded set of $\R^n$.

We denote by $\R^{n\times n}$ the set of $n\times n$ matrices, by $\Gamma_{+}$ its subset of symmetric positive semidefinite matrices, by $SO(n) \subset \R^{n\times n}$ the orthogonal matrices with determinant $1$, and by $I$ the identity matrix.

Given $A_i \in \R^{n\times n}$, the measure $\delta_{A_i}$ is the Dirac delta at $A_i$.
The barycenter of the probability measure $\nu=\sum_{i=1}^{N}\alpha_{i}\delta_{A_{i}}$ is
$\overline{\nu}=\sum_{i=1}^{N}\alpha_{i}A_{i}$. 

Given $A\in \R^{n\times n}$, let $\sigma_1(A)\leq\cdots\leq\sigma_n (A)$ denote its singular values. If the matrix $A$ is clear from the context, we will just indicate their singular values as $\sigma_1,\ldots,\sigma_n.$
In fact, in this paper we will always deal with $A \in \Gamma_+$, so its eigenvalues coincide with its singular values.
Its components are written $A_{\alpha, \beta}$ for $\alpha , \beta \in \{ 1, \ldots, n \}$.
Its operator norm is denoted by $|A|$, which coincides with $\sigma_n (A)$.
The norm of a $v \in \R^n$ is also denoted by $|v|$.

Given $a_1, \ldots, a_n \in \R$ the matrix $\diag (a_1, \ldots, a_n) \in \R^{n\times n}$ is the diagonal matrix with diagonal entries $a_1, \ldots, a_n$.

We will use the symbol $\lesssim$ when there exists a constant depending only on $n$, $m_1$ and $m_2$ such that the left hand side is less than or equal to the constant times the right hand side. Sometimes, the left hand side could be negative.

Given a set $E \subset \R^n$, we denote its characteristic function by $\chi_E$. We write $\# E$ for the number of elements of $E$. When it is measurable, its Lebesgue measure is denoted by $|E|$ and we use $\mathcal{H}^{m}(E)$ for its Hausdorff measure of dimension $m$.

Given $a \in \R$, its integer part is denoted by $\floor{a}$ and we denote by $\ceil{a}$ its ceiling function.

Given $E \subset \R^n$, $\alpha \in (0,1]$ and a function $f : E \to \R^n$, we denote the H\"older seminorm, supremum norm and H\"older norm, respectively, as
\begin{align*}
 & \left| f \right|_{C^{\alpha} (E)} := \sup_{\substack{x_1, x_2 \in E \\ x_1 \neq x_2}} \frac{\left| f(x_2) - f(x_1)\right|}{\left| x_2 - x_1 \right|^{\alpha}} , \qquad \left\| f \right\|_{L^{\infty} (E)} = \sup_{x \in E} \left| f(x) \right| , \\
 & \left\| f \right\|_{C^{\alpha} (E)} :=  \left| f \right|_{C^{\alpha} (E)} + \left\| f \right\|_{L^{\infty} (E)} .
\end{align*}
We will write $f \in C^{\alpha} (E, \R^n)$ when $\left\| f \right\|_{C^{\alpha} (E)} < \infty$. 
Note that, if $f$ is continuous up to the boundary, the above norms and seminorms in $E$ coincide with those in $\overline{E}$.
In particular, we will identify $C^{\alpha} (E, \R^n)$ with $C^{\alpha} (\overline{E}, \R^n)$, the set of H\"older functions of exponent $\alpha$.
Of course, if $\alpha=1$, they are Lipschitz. 

The identity function is denoted by $\id$ and the Sobolev space from $\Omega$ to $\R^{n}$ is denoted, alternatively, by $W^{1,p}$, $W^{1,p}(\Omega)$ or $W^{1,p}(\Omega,\R^{n})$.

Given $f:A\to\R^{n}$, where $A$ is a subset of an $m$-dimensional affine space of $\R^{n}$, we say that it satisfies the $m$-dimensional Luzin condition, also known as condition $N$, if for every $E\subset A$ with $\mathcal{H}^{m}(E)=0$, then $\mathcal{H}^{m}(f(E))=0$.

For $f\in W^{1,1}(\Omega,\R^{n})$, we denote the Jacobian of $f$ by $Jf$ and, for $k\in\{1,\ldots,n\}$, the $k$-dimensional Jacobian of $f$ by $J_{k}f$, i.e.,  $J_{k}f$ is equal to $\left(\sum_{M}|M|^{2}\right)^{\frac{1}{2}}$, where the sum runs over all the minors of $Df$ of order $k$.

We will say that a continuous mapping $f:\overline{\Omega}\to \R^n$ is \emph{piecewise affine} if there exists a countable family $\{\Omega_i\}_{i\in \N}$ of pairwise disjoint open subsets of $\Omega$ such that $f|_{\Omega_i}$ is affine for all $i\in\N$, and 
\[\left|\Omega\setminus \bigcup_{i\in \N}\Omega_i\right|=0.\]
Note that $\{\Omega_{i}\}_{i\in\N}$ need not be locally finite. Given $S\subset \R^{n\times n}$ a set of invertible matrices we denote by $S^{-1}$ the set
\[\{A^{-1}: A\in S\}.\]
\section{Limitations to the integrability of $f$ and $f^{-1}$}\label{sect: sharp}
In this section we prove Theorems \ref{theorem integrability of f and rank of the inverse}, \ref{theorem sum of the integrabilities} and \ref{theorem integrabilities and rank}.  The proof of Theorem \ref{theorem integrability of f and rank of the inverse} follows that in \cite[Th. 4]{HeMoPaSb09}.
\begin{proof}[Proof of Theorem \ref{theorem integrability of f and rank of the inverse}]
Suppose that $\rank(Df)\geq n-m$ in a set $A\subset E$ of positive measure. Then $|J_{n-m}f|>0$ on $A$, and without loss of generality we can assume that $Jf=0$ on $A$ and that $f$ is Lipschitz on $A$, see \cite[Section 6.6, Th.3]{EvGa92}.

For each $I\subset\{1,\ldots,n\}$ with $|I|=n-m$ let $\pi_{I}:\R^{n}\to\R^{n-m}$ be the projection that sends $x=(x_1,\ldots,x_n)$ to $(x_{i_1},\ldots,x_{i_{n-m}})$ , where $I=\{i_1,\ldots,i_{n-m}\}$ and $i_{1}<\cdots<i_{n-m}$. Define $h_I=\pi_I\circ f$ and set $P(z)=\pi_I^{-1}(z)\cap f(\Omega)$ for $z\in\R^{n-m}$. Since $|J_{n-m}f|>0$ in $A$, there exists $I\subset \{1,\ldots,n\}$ such that $|I|=n-m$ and $|J_{n-m}h_I|>0$ on a subset of $A$ of positive measure, still called $A$. Since $f$ is Lipschitz on $A$, we can use the area formula (see, e.g. \cite{EvGa92}, \cite{Federer69}, \cite{AmFuPa00}) to conclude that
\[|f(A)|=0.\]
Using the coarea formula we get 
\[0<\int_{A}|J_{n-m}h_I|=\int_{\R^{n-m}} \mathcal{H}^{m}\left( \{x\in A: h_I(x)=z\} \right)dz.\]
Hence, for a set $F\subset\R^{n-m}$ of positive measure we have that $\forall z\in F$,
\begin{equation}\label{proof sharp >0}
\mathcal{H}^{m}\left(f^{-1}\left(f(A)\cap P(z)\right) \right)=\mathcal{H}^{m}\left( \left\{ x\in A: h_I(x)=z \right\}\right)>0.
\end{equation}
In the equality we have used that $f$ is injective. On the other hand, since $|f(A)|=0$ it follows that for $\mathcal{H}^{n-m}$-almost every $z\in\R^{n-m}$ we get $\mathcal{H}^{m}\left(f(A)\cap P(z)\right)=0$. 

As $f^{-1}\in W^{1,q}$, we have that $f^{-1}\in W^{1,q}(P(z))$ for $\mathcal{H}^{n-m}$-almost every $z\in\R^{n-m}$, so under either option \emph{\ref{optionA})} or \emph{\ref{optionB})}, $f^{-1}|_{P(z)}$ satisfies the $m$-dimensional Luzin (N) condition. The proof under \emph{\ref{optionA})} is due to \cite[Theorem 1.1]{KoMaZu15} (with $\lambda=0$ in the notation there), while the proof under \emph{\ref{optionB})} is classical \cite{MaMi73}. Therefore, for $\mathcal{H}^{n-m}$-almost every $z\in\R^{n-m}$, we obtain
\[\mathcal{H}^{m}\left(f^{-1}\left(f(A)\cap P(z)\right) \right)=0,\]
and we have a contradiction with (\ref{proof sharp >0}).

\end{proof}
Now using Theorem \ref{theorem integrability of f and rank of the inverse} we are able to prove Theorem \ref{theorem sum of the integrabilities}.
\begin{proof}[Proof of Theorem \ref{theorem sum of the integrabilities}]
Define $m_p=\ceil{p}-1$ and $m_q=\ceil{q}-1$. Observe that $m_p< p\leq m_p+1$ and $m_q< q\leq m_q+1$. Then, using Theorem \ref{theorem integrability of f and rank of the inverse}, we have that $\rank Df\leq n-m_q-1$ almost everywhere on $\Omega$. On the other hand, using \cite[Th. 12]{FaMoOl16}, we get that $\rank Df\geq m_p$ in a subset of $\Omega$ of positive measure; we also obtain $\rank Df^{-1}\geq m_q$. Therefore, we have
\[m_p\leq n-m_q-1.\]
Hence
\begin{equation}\label{proof sum of int p+q<n+1}
p+q\leq m_p+m_q+2\leq n+1.
\end{equation}
If $p+q=n+1$ then $p=m_p+1$ and $q=m_q+1$. If, in addition, $f^{-1}$ or $f$ were H\"older continuous then, by Theorem \ref{theorem integrability of f and rank of the inverse}, $\rank( Df)\leq n-q-1$ or $\rank (Df^{-1})\leq n-p-1$, but this contradicts $\rank(Df)\geq m_p=p-1$ in the first case and $\rank(Df^{-1})\geq m_q=q-1$ in the second case.
\end{proof}
\begin{proof}[Proof of Theorem \ref{theorem integrabilities and rank}]
First, we observe that thanks to \cite[Th. 12]{FaMoOl16} we obtain $p\leq m_1+1$ and $q\leq m_2+1$.

On the other hand, denote $m=\ceil{q}-1$, then $m< q\leq m+1$, and using Theorem \ref{theorem integrability of f and rank of the inverse} with $E=A$ we obtain $\rank(Df)\leq n-m-1\leq n-q$ almost everywhere on $A$, so $q\leq n-m_1$. In the same way, we get $p\leq n-m_2$. The theorem follows.
\end{proof}
\section{Construction of the laminate in dimension three}\label{sect: sketch proof}

The next definition introduces the concept of laminate of finite order \cite{Dacorogna89, Pedregal93, MuSv03, AsFaSz08}. 

\begin{defn}\label{de:laminate}

The family $\mathcal{L}(\R^{n\times n})$ of laminates of finite order is the smallest family of probability measures in $\R^{n\times n}$ with the properties:

\begin{enumerate}[i)]
\item $\mathcal{L}(\R^{n\times n})$ contains all the Dirac masses.
\item\label{def lam split proc} If $\sum_{i=1}^{N}\lambda_i\delta_{A_i}\in \mathcal{L}(\R^{n\times n})$ and $A_N=\lambda B+(1-\lambda)C$, where $\lambda\in [0,1]$ and $\rank \,(B-C)=1$, then the probability measure 
\[\sum_{i=1}^{N-1}\lambda_i\delta_{A_i}+\lambda_N(\lambda \delta_B+(1-\lambda)\delta_C)\]
is also in $\mathcal{L}(\R^{n\times n})$. 
\end{enumerate}
\end{defn}

Since in this work we will only use laminates of finite order, for simplicity they will be just called \emph{laminates}. Note that they are a convex combination of Dirac masses.

The next lemma gives us a characterization of the laminates of finite order.
\begin{lem}\label{lemma caract of prelaminates}
For every laminate of finite order $\nu=\sum_{i=1}^{N}\lambda_i \delta_{A_{i}}$ there exists a family $\{\nu_j\}_{j=1}^{N}$ of laminates of finite order, such that
\begin{itemize}
\item $\nu_1=\delta_{\overline{\nu}}$,
\item for $j\in\{1,\ldots,N-1\}$, $\nu_{j+1}$ is obtained from $\nu_j$ using once \ref{def lam split proc}) of Definition \ref{de:laminate},
\item $\nu_{N}=\nu$.
\end{itemize}
\end{lem}
Lemma \ref{lemma caract of prelaminates} has an elementary proof. It is illustrated in the following example. We have the splits

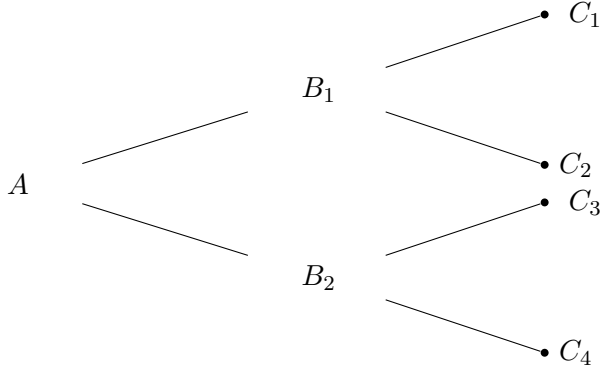
\begin{tikzpicture}[grow=right, sloped]
\node[bag] {$A$}
  child {
        node[bag] {$B_2$}        
        child {
                node[end, label=right:
                    {$C_4$}] {}
            }
            child {
                node[end, label=right:
                    { $C_3$}] {}
            }
            }
    child {
        node[bag] {$B_1$}        
        child {
                node[end, label=right:
                    {$C_2$}] {}
            }
            child {
                node[end, label=right:
                    { $C_1$}] {}
            }
    };
\end{tikzpicture}

and $\nu=\sum_{i=1}^{4}\lambda_i\delta_{C_i}$. Then the laminates $\{\nu_{i}\}_{i=1}^{4}$ can be chosen to be the following:
\[\nu_4=\nu,\]
\[\nu_3=(\lambda_1+\lambda_2)\delta_{B_1}+\lambda_3\delta_{C_3}+\lambda_{4}\delta_{C_4},\]
\[\nu_2=(\lambda_1+\lambda_2)\delta_{B_1}+(\lambda_3+\lambda_{4})\delta_{B_2},\]
\[\nu_1=\delta_{A}.\]
From Lemma \ref{lemma caract of prelaminates} we obtain the following corollary that will be used throughout the paper to prove that some measures are laminates.
\begin{cor}\label{cor: laminate}
Let $\nu=\sum_{i=1}^{N}\lambda_i\delta_{A_{i}}$ and $\{\nu_{A_i}\}_{i=1}^{N}$ be laminates of finite order such that
\[\overline{\nu}_{A_i}=A_i\text{ with }A_i \text{ being all different}.\]
Then, the probability measure
\[\nu'=\sum_{i=1}^{N}\lambda_i\nu_{A_{i}}=\nu+\sum_{i=1}^{N}\nu(A_i)[\nu_{A_i}-\delta_{A_i}]\]
is also a laminate of finite order.
\end{cor}
\begin{defn}\label{de:inverse laminate}
We define the inverse laminate of a laminate $\nu=\sum_{i=1}^{N}\lambda_i \delta_{A_i}$ supported in the set of positive definite matrices as
\[\nu^{-1}=\frac{1}{\det(\overline{\nu})}\sum_{i=1}^{N}\lambda_i \det(A_i)\delta_{A_i^{-1}}.\]
\end{defn}
This definition, which seems to be new, arises from the fact that if $\sum_{i=1}^{N}\lambda_i \delta_{A_i}$ is a laminate supported in the set of positive definite matrices and $f$ is a piecewise affine Sobolev homeomorphism satisfying 
\[|E_i|=\lambda_i\text{ and }f_{i}|_{\partial E_{i}}=A_{i}\qquad\text{for }i\in\{1,\ldots,N\},\]
for 
\[E_{i}=\{x\in \Omega: |Df(x)-A_i|<\delta\},\]
 and some $\delta>0$, then we get 
\[|f(E_{i})|=\int_{E_{i}}\det Df(x)dx=\int_{E_{i}}\det A_{i}dx=\lambda_i\det A_i\qquad\text{for }i\in\{1,\ldots,N\},\]
and hence, there exists $\delta'>0$ such that $f^{-1}$ satisfies
\[|\{y\in f(\Omega): |Df^{-1}(y)-A_{i}^{-1}|<\delta'\}|=|f(E_{i})|=\lambda_{i}\det A_{i}.\]
Although we will not use it, $\nu^{-1}$ is also a laminate; this can be shown using the fact that $\det$ is rank-one linear.

In this section we construct the sequence of laminates $ \nu_{j}$ of finite order that is behind the proof of Theorem \ref{theorem existence} in the case $n=3$ and $m_1=m_2=1$. The actual proof will consist in approximating $\nu_{k}$ with laminates of finite order supported in the set of positive definite matrices, then use Proposition \ref{prop: laminate-bisobolev} to obtain homeomorphisms that are close to the approximate laminates and such that their inverse are close to the inverse of that laminates, then paste the obtained homeomorphisms to construct a homeomorphism in the whole domain and, finally, a limit passage will yield the homeomorphism $f$ of Theorem \ref{theorem existence}.

Although this section is not strictly necessary for the proof of Theorem \ref{theorem existence}, it will help the reader to follow the construction of Section \ref{sect: proof of existence}.

To define $\nu_j$ we need to define the following sets. For $i,k\in\N$ let
\begin{equation}\label{definicion Aki n=3}
A_{k}^{i}=\left\{A\in \Gamma_{+}: \sigma_j(A)=k^{-1} \text{ for }j\in\{1,2\}\text{ and }\sigma_3(A)\in\{i-1,i\}\setminus\{0\}\right\},
\end{equation}
\begin{equation}\label{definicion Bki n=3}
B_{k}^{i}=\left\{A\in \Gamma_{+}: \sigma_1(A)\in\{k^{-1},(k-1)^{-1}\}\setminus\{\infty\}\text{ and }\sigma_j(A)=i\text{ for }j\in\{2,3\}\right\}\setminus\{I\}.
\end{equation}
Observe that a matrix $A\in\Gamma_{+}\setminus\{I\}$ belongs to $A_{k}^{i}$ if and only if $A^{-1}$ is in $B_{i}^{k}$.

The laminates $\nu_j$ that we will construct will satisfy the following:
\begin{enumerate}[(a)]
\item\label{sk ident} $\overline{\nu}_j=I$,
\item\label{sk sup} $\supp(\nu_j)\subset \bigcup_{i=1}^{j}A_j^{i}\cup B_{i}^{j}$.

\item\label{sk int of f}For all $\ve>0$ there exists a bounded family of constants $\{C_{j}\}_{j\in\N}$, such that, for all $j\in\N$,
\[\nu_j(A_j^{i})\leq C_j i^{-3+\ve} \text{ and }\nu_j(B_i^{j})\leq C_j i^{-2+\ve}j^{-2}.\]
\end{enumerate}
When we approximate these laminates by functions, $f_j$, and then pass to the limit, we obtain a bi-Sobolev homeomorphism $f$ that is the identity on the border due to (\ref{sk ident}); by (\ref{sk sup}) and (\ref{sk int of f}) we get for $f_j$
\[Df_j\in \bigcup_{i=1}^{j}\left( A_{j}^{i}\cup B_{i}^{j}\right)+B(0,r_j)\text{ with }r_j\to 0,\]
\[|\{x: Df_{j}(x)\in \bigcup_{i=1}^{j} B_{i}^{j}+B(0, r_j)\}|\to 0,\]
and, for the inverse, using that $A\in A_{j}^{i}$ if and only if $A^{-1}\in B_{i}^{j}$, we obtain
\[Df_{j}^{-1}\in \bigcup_{i=1}^{j}\left( A_{j}^{i}\cup B_{i}^{j}\right)+B(0,r_j')\text{ with }r_j'\to 0,\]
\[|\{y: Df_{j}^{-1}(y)\in \bigcup_{i=1}^{j} B_{i}^{j}+B(0, r_j')\}|\to 0.\]
So, the ranks of $Df$ and $Df^{-1}$ are equal to $1$ almost everywhere. Moreover, with (\ref{sk int of f}) we obtain the following.
Let $t,\ve>0$, and pick $j\in\N$, $j>t$ and big enough; then
\[\nu_j\left(\{A\in\R^{3 \times 3}: |A|>t\}\right)\lesssim \sum_{i=\ceil{t}}^{j}i^{-3+\ve}+j^{-2}\sum_{i=1}^{j}i^{-2+\ve}\lesssim t^{-2+\ve},\]
\begin{align*}
\nu_j^{-1}\left(\{A\in\R^{3 \times 3}: |A|>t\}\right)&\lesssim \sum_{i=1}^{j}i^{-3+\ve}\sup_{M\in A_{j}^{i}}\{\det M\}+j^{-2}\sum_{i=\ceil{t}}^{j}i^{-2+\ve}\sup_{M\in B_{i}^{j}}\{\det M\}\\
&\leq j^{-2}\sum_{i=1}^{j}i^{-2+\ve}+2\sum_{i=\ceil{t}}^{j}i^{-3+\ve} \lesssim t^{-2+\ve}.
\end{align*}
Therefore, this gives us that $Df, Df^{-1}\in \bigcap_{p<2}W^{1,p}$, respectively.

The laminates $\nu_j$ are defined inductively as follows. We begin with $\nu_1=I$. It is clear that $\nu_1$ satisfies (\ref{sk ident}), (\ref{sk sup}) and (\ref{sk int of f}). Now, let
\begin{equation}\label{sk nuj}
\nu_{j}=\sum_{k=1}^{N}\lambda_{k}\delta_{A_k}\in\mathcal{L}\left(\R^{n\times n}\right),
\end{equation}
with $A_k\in \bigcup_{i=1}^{j}( A_{j}^{i}\cup B_{i}^{j})$, all different. For each $A\in\supp(\nu_{j})$ we are going to construct a laminate $\nu_{A}$, whose support is in $\bigcup_{i=1}^{j+1}( A_{j+1}^{i}\cup B_{i}^{j+1})$. To do that, we need the following two lemmas.
\begin{lem}\label{lemma sk laminado Aki}
%En el paper de verdad mejor hacer que este cerca de Bki y hacer el mismo laminado
Let $A\in A_{k}^{i}$. Then there exists a laminate of finite order $\nu$ such that 
\begin{itemize}
\item $\overline{\nu}=A$,
\item $\supp (\nu)\subset A_{k+1}^{i}\cup B_{k+1}^{i}$,
\item $\nu(B_{k+1}^{i})\lesssim \left(k^{2}i\right)^{-1}$.
\end{itemize}
\end{lem}
\begin{lem}\label{lemma sk laminado Bki}
Let $A\in B_{k}^{i}$. Then there exists a laminate of finite order $\nu$ such that 
\begin{itemize}
\item $\overline{\nu}=A$,
\item $\supp (\nu)\subset A_{k}^{i+1}\cup B_{k}^{i+1}$,
\item $\nu(A_{k}^{i+1})\lesssim i^{-1}$,
\item $\nu(B_{k}^{i+1})- \left(\frac{i}{i+1}\right)^{2}\lesssim (ik)^{-2}$.
\end{itemize} 
\end{lem}
We will only prove Lemma \ref{lemma sk laminado Aki}, the proof of Lemma \ref{lemma sk laminado Bki} being analogous.
\begin{proof}[Proof of Lemma \ref{lemma sk laminado Aki}]

Given $A\in A_{k}^{i}$, without loss of generality, we can assume that $A=\diag(k^{-1},k^{-1},\sigma)$, with $\sigma\in \{\max\{i-1,1\},i\}$. 

Now, we denote $\mu=\frac{i-k^{-1}}{i-(k+1)^{-1}}$; observe that $\mu$ satisfies
\begin{equation}\label{pro Aki mu}
0<\mu\leq 1\qquad \text{and }1-\mu\leq (k^{2}i)^{-1}.
\end{equation}
Using that
\begin{equation}\label{pro Aki kmu}
k^{-1}=\mu (k+1)^{-1}+(1-\mu)i,
\end{equation}
we split $A$ in the following way:

\begin{tikzpicture}[grow=right, sloped]
\node[bag] {$(k^{-1},k^{-1},\sigma)$}
    child {
        node[end, label=right: {$(i,k^{-1},\sigma).$}]{}        
            edge from parent 
            node[below]  {$1-\mu$}
    }
    child {
        node[bag] {$((k+1)^{-1},k^{-1},\sigma)$}        
        child {
                node[end, label=right:
                    {$((k+1)^{-1},i,\sigma)$}] {}
                edge from parent
                node[below]  {$1-\mu$}
            }
            child {
                node[end, label=right:
                    { $((k+1)^{-1},(k+1)^{-1},\sigma)\in A_{k+1}^{i}$}] {}
                edge from parent
                node[above] {$\mu$}
            }
        edge from parent         
            node[above] {$\mu$}
    };
\end{tikzpicture}\\
Therefore
\begin{align*}
A=&\mu^{2}\diag((k+1)^{-1},(k+1)^{-1},\sigma)+\mu(1-\mu)\diag((k+1)^{-1},i,\sigma)\\
&+(1-\mu)\diag(i,k^{-1},\sigma).
\end{align*}
If $\sigma=i$ we define 
\[\nu=\mu^{2}\delta_{\diag((k+1)^{-1},(k+1)^{-1},i)}+\mu(1-\mu)\delta_{\diag((k+1)^{-1},i,i)}
+(1-\mu)\delta_{\diag(i,k^{-1},i)},\]
which, clearly, is a laminate supported in $A_{k+1}^{i}\cup B_{k+1}^{i}$, and by (\ref{pro Aki mu}), we have that
\[\nu(B_{k+1}^{i})\lesssim (k^{2}i)^{-1}.\]
If, on the contrary, $\sigma=i-1$, we define $\lambda=\frac{i-1-(k+1)^{-1}}{i-(k+1)^{-1}}$ and using (\ref{pro Aki kmu}) and that 
\[i-1=\lambda i+(1-\lambda)(k+1)^{-1},\]
we do the following splits:

\begin{tikzpicture}[grow=right, sloped]
\node[bag] {$((k+1)^{-1},i,i-1)$}
    child {
        node[end, label=right: {$((k+1)^{-1},i,(k+1)^{-1})\in A_{k+1}^{i},$}]{}        
            edge from parent 
            node[below]  {$1-\lambda$}
    }
    child {
        node[end, label=right: {$((k+1)^{-1},i,i)\in B_{k+1}^{i}$}]{}        
        edge from parent         
            node[above] {$\lambda$}
    };
\end{tikzpicture}

\begin{tikzpicture}[grow=right, sloped]
\node[bag] {$(i,k^{-1},i-1)$}
    child {
        node[bag] {$(i,k^{-1},(k+1)^{-1})$}    
        child {
                node[end, label=right:
                    {$(i,i,(k+1)^{-1})\in B_{k+1}^{i}.$}] {}
                edge from parent
                node[below]  {$1-\mu$}
            }
            child {
                node[end, label=right:
                    { $(i,(k+1)^{-1},(k+1)^{-1})\in A_{k+1}^{i}$}] {}
                edge from parent
                node[above] {$\mu$}
            }
            edge from parent 
            node[below]  {$1-\lambda$}
    }
    child {
        node[end, label=right:
            {$(i,k^{-1},i)\in B_{k+1}^{i}$}]  {}        
        edge from parent         
            node[above] {$\lambda$}
    };
\end{tikzpicture}\\
Hence
\begin{align*}
A=&\mu^{2}\diag((k+1)^{-1},(k+1)^{-1},i-1)+\mu(1-\mu)\lambda\diag((k+1)^{-1},i,i)\\
&+\mu(1-\mu)(1-\lambda)\diag((k+1)^{-1},i,(k+1)^{-1})\\
&+(1-\mu)\lambda\diag(i,k^{-1},i)+(1-\mu)(1-\lambda)\mu\diag(i,(k+1)^{-1},(k+1)^{-1})\\
&+(1-\mu)^{2}(1-\lambda)\diag(i,i,(k+1)^{-1}),
\end{align*}
and we define
\begin{align*}
\nu=&\mu^{2}\delta_{\diag((k+1)^{-1},(k+1)^{-1},i-1)}+\mu(1-\mu)\lambda\delta_{\diag((k+1)^{-1},i,i)}\\
&+\mu(1-\mu)(1-\lambda)\delta_{\diag((k+1)^{-1},i,(k+1)^{-1})}\\
&+(1-\mu)\lambda\delta_{\diag(i,k^{-1},i)}+(1-\mu)(1-\lambda)\mu\delta_{\diag(i,(k+1)^{-1},(k+1)^{-1})}\\
&+(1-\mu)^{2}(1-\lambda)\delta_{\diag(i,i,(k+1)^{-1})},
\end{align*}
which is a laminate supported in $A_{k+1}^{i}\cup B_{k+1}^{i}$. Using (\ref{pro Aki mu}) we obtain
\[\nu (B_{k+1}^{i})\leq 1-\mu^{2}\lesssim 1-\mu\leq (k^{2}i)^{-1},\]
and the proof is complete.
\end{proof}
Now, we can prove the next two lemmas that will give us the laminate $\nu_A$.

\begin{lem}\label{lemma sk laminado Aji}
Let $i\in\N$, $i\leq j$ and $A\in A_j^{i}$. Then, there exists a laminate $\nu$ such that 
\begin{itemize}
\item $\overline{\nu}=A$,
\item $\supp (\nu)\subset \bigcup_{b=0}^{j-i+1}A_{j+1}^{i+b}\cup B_{j+1}^{j+1}$,
\item $\nu(A_{j+1}^{i+b})\lesssim j^{-2}  \frac{i}{(i+b)^{3} }$, for $b\in\{1,\ldots,j-i+1\}$,
\item $\nu(B_{j+1}^{j+1})\lesssim j^{-2}i(j+1)^{-2}$.
\end{itemize}
\end{lem}

\begin{lem}\label{lemma sk laminado Bkj}
Let $i\in\N$, $i\leq j$ and $A\in B_i^{j}$. Then, there exists a laminate $\nu$ such that 
\begin{itemize}
\item $\overline{\nu}=A$,
\item $\supp (\nu)\subset A_{j+1}^{j+1}\cup\bigcup_{b=0}^{j-i+1} B_{i+b}^{j+1}$,
\item $\nu(A_{j+1}^{j+1})\lesssim  j^{-1}$, 
\item $\nu(B_{i}^{j+1})- \left(\frac{j}{j+1}\right)^{2}\lesssim j^{-2}i^{-2}$,
\item $\nu(B_{i+b}^{j+1})\lesssim \left((i+b-1)(j+1)\right)^{-2} $, for $b\in\{1,\ldots,j-i+1\}$.
\end{itemize}
\end{lem}

As before, we will only prove Lemma \ref{lemma sk laminado Aji} since the proof of Lemma \ref{lemma sk laminado Bkj} can be obtained in the same form.
\begin{proof}[Proof of Lemma \ref{lemma sk laminado Aji}]
It is enough to construct a family of laminates $\{\nu_{\ell}'\}_{\ell=1}^{j-i+2}$ such that
\begin{enumerate}[i)]
\item\label{sk Aji id} $\overline{\nu}_{\ell}'=A$,
\item\label{sk Aji sup}  $\supp (\nu_{\ell}')\subset \bigcup_{b=0}^{\ell-1}A_{j+1}^{i+b}\cup B_{j+1}^{i+\ell-1}$,
\item\label{sk Aji bound Aj+1i+b}  $\nu_{\ell}'(A_{j+1}^{i+b})\lesssim j^{-2}  \frac{i}{(i+b)^{3} }$, for $b\in\{1,\ldots,\ell-1\}$,
\item\label{sk Aji bound Bj+1i+l-1}  $\nu_{\ell}'(B_{j+1}^{i+\ell-1})\lesssim j^{-2}i(i+\ell-1)^{-2}\left(1+j^{-2}8C\sum_{k=1}^{\ell-1}(i+k-1)^{-2}\right)$,
\end{enumerate}
and define $\nu=\nu_{j-i+2}'$. The constant $C$ is bigger than those that appear in Lemma \ref{lemma sk laminado Bki}. Let $\nu_1'$ be the laminate of Lemma \ref{lemma sk laminado Aki}; then, $\nu_1'$ satisfies all the conditions. Now, for $1\leq \ell\leq j-i+1$, given $\nu_{\ell}'=\sum_{a=1}^{N_{\ell}}\lambda_{a}\delta_{B_{a}}$ with all the $B_{a}$ different, define $\nu_{B_{a}}$ as the laminate of Lemma \ref{lemma sk laminado Bki} if $B_{a}\in B_{j+1}^{i+\ell-1}$ and as $\delta_{B_{a}}$ otherwise. Set 
\[\nu_{\ell+1}'=\sum_{a=1}^{N_{\ell}}\lambda_{a}\nu_{B_{a}}.\]
It is immediate that $\nu_{\ell+1}'$ satisfies \ref{sk Aji id}), \ref{sk Aji sup}), and \ref{sk Aji bound Aj+1i+b}) for $b\in\{1,\ldots,\ell-1\}$. Thanks to Corollary \ref{cor: laminate}, it is a laminate. Hence, we only have to bound $\nu_{\ell+1}'(A_{j+1}^{i+\ell})$ and $\nu_{\ell+1}'(B_{j+1}^{i+\ell})$. We have
\begin{align*}
\nu_{\ell+1}'(A_{j+1}^{i+\ell})&\lesssim \nu_{\ell}'(B_{j+1}^{i+\ell-1})(i+\ell-1)^{-1}\lesssim j^{-2}i(i+\ell-1)^{-3}\lesssim j^{-2}  \frac{i}{(i+\ell)^{3} }.
\end{align*}
For $j$ big enough, $j^{-2}8C\sum_{k=1}^{\infty}(i+k-1)^{-2}\leq 1$, so
\begin{align*}
&1+4C((i+\ell-1)(j+1))^{-2}+\left( 1+4C((i+\ell-1)(j+1))^{-2}\right) j^{-2}8C\sum_{k=1}^{\ell-1}(i+k-1)^{-2}\\
&\leq 1+j^{-2}8C\sum_{k=1}^{\ell}(i+k-1)^{-2}.
\end{align*}
Hence,
\begin{align*}
\nu_{\ell+1}'(B_{j+1}^{i+\ell})&\leq \nu_{\ell}'(B_{j+1}^{i+\ell-1})\left(\left(\frac{i+\ell-1}{i+\ell}\right)^{2}+C((i+\ell-1)(j+1))^{-2}\right)\\
&\lesssim j^{-2}i(i+\ell)^{-2}\Biggl(1+4C((i+\ell-1)(j+1))^{-2}\bigr.\\
&\quad\left.+\left( 1+4C((i+\ell-1)(j+1))^{-2}\right) j^{-2}8C\sum_{k=1}^{\ell-1}(i+k-1)^{-2}\right)\\
&\leq j^{-2}i(i+\ell)^{-2}\left(1+j^{-2}8C\sum_{k=1}^{\ell}(i+k-1)^{-2}\right).
\end{align*}
\end{proof}
From $\nu_j$ as in (\ref{sk nuj}), we construct $\nu_{j+1}$ as follows. For each $A\in\supp (\nu_j)$ we define $\nu_{A}$ as the laminate of Lemma \ref{lemma sk laminado Aji} if $A\in \bigcup_{i=1}^{j}A_{j}^{i}$, and as the laminate of Lemma \ref{lemma sk laminado Bkj} if $A\in \bigcup_{i=1}^{j}B_{i}^{j}$. We define $\nu_{j+1}$ as
\[\nu_{j+1}=\sum_{k=1}^{N}\lambda_k \nu_{A_k}.\]
It is clear that $\nu_{j+1}$ satisfies (\ref{sk ident}) and (\ref{sk sup}). Pick $\ve>0$ small enough to satisfy $\sum_{a=1}^{\infty}a^{-2+\ve}\leq 2$. We will prove by induction on $j$ that $\nu_{j+1}$ satisfies (\ref{sk int of f}). Let $C_0>0$ be a constant bigger than the constants appearing in Lemmas \ref{lemma sk laminado Aji} and \ref{lemma sk laminado Bkj}. Set $j_{\ve}\in\N$ such that $j_{\ve}^{\ve}>20C_0$. Define 
\[C_{j}=\max_{i\in\{1,\ldots,j\}}\left\{\nu_{j}(A_{j}^{i})i^{3-\ve},\nu_{j}(B_{i}^{j})i^{2-\ve}j^{2}\right\}\quad\text{for }j\leq j_{\ve},\]
and 
\[C_{j+1}=C_{j}(1+12C_0j^{-2})\quad \text{for }j\geq j_{\ve}.\]
It is clear that $\sup_{j\in\N}C_j<\infty,$ and we have (\ref{sk int of f}) for $j\leq j_{\ve}$.

Then for $j\geq j_{\ve}$ and for each $i\in\{1,\ldots,j\}$ we know that the matrices in $A_{j+1}^{i}$ can only come from $\bigcup_{a=1}^{i}A_{j}^{a}$, and the matrices in $B_{i}^{j+1}$ can only come from $\bigcup_{a=1}^{i}B_{a}^{j}$; therefore we obtain
\begin{align*}
\nu_{j+1}(A_{j+1}^{i})&\leq \nu_{j}(A_{j}^{i})+\sum_{a=1}^{i-1}\nu_{j}(A_{j}^{a})C_0 j^{-2}\frac{a}{i^{3}}\leq C_j \left(i^{-3+\ve}+\sum_{a=1}^{i-1}C_0 j^{-2} a^{-2+\ve} i^{-3} \right)\\
&\leq C_j (1+2C_0 j^{-2}) i^{-3+\ve}\leq C_{j+1}i^{-3+\ve},
\end{align*}
\begin{align*}
\nu_{j+1}(B_{i}^{j+1})&\leq \nu_{j}(B_{i}^{j})\left(\left(\frac{j}{j+1}\right)^{2}+C_0 j^{-2}i^{-2}\right)+\sum_{a=1}^{i-1}\nu_{j}(B_{a}^{j})C_0 ((i-1)(j+1))^{-2}\\
&\leq C_j \left(i^{-2+\ve}(j+1)^{-2}+C_0 i^{-4+\ve} j^{-4}+4C_0  i^{-2}j^{-2}(j+1)^{-2}\sum_{a=1}^{i-1}a^{-2+\ve}\right)\\
&\leq C_j i^{-2+\ve}(j+1)^{-2}(1+12C_0  j^{-2})= C_{j+1}i^{-2+\ve}(j+1)^{-2}.
\end{align*}
For $i=j+1$, since the matrices in $A_{j+1}^{j+1}\cup B^{j+1}_{j+1}$ can come from any matrix in the support of $\nu_j$ we get
\begin{align*}
\nu_{j+1}(A_{j+1}^{j+1})&\leq \sum_{a=1}^{j}\left[\nu_{j}(A_{j}^{a})C_0 j^{-2}\frac{a}{(j+1)^{3}}+\nu_{j}(B_{a}^{j})C_0 j^{-1}\right]\\
&\leq C_jC_0\sum_{a=1}^{i-1}\left[ j^{-2} a^{-2+\ve} (j+1)^{-3} +a^{-2+\ve}j^{-3}\right]\\
&\leq C_j  20 C_0(j+1)^{-3}\leq C_{j+1}(j+1)^{-3+\ve},
\end{align*}
\begin{align*}
\nu_{j+1}(B_{j+1}^{j+1})&\leq \sum_{a=1}^{j}\left[\nu_{j}(A_{j}^{a})C_0 j^{-2}a(j+1)^{-2}+\nu_{j}(B_{a}^{j})C_0 j^{-2}(j+1)^{-2}\right]\\
&\leq C_j\left(C_0\sum_{a=1}^{i-1}\left[ j^{-2} a^{-2+\ve} (j+1)^{-2} +a^{-2+\ve}j^{-4}(j+1)^{-2}\right]\right)\\
&\leq C_j 20 C_0 (j+1)^{-4}\leq C_{j+1}(j+1)^{-4+\ve},
\end{align*}
and the proof of (\ref{sk ident})--(\ref{sk int of f}) is completed.
\section{Proof of Theorem \ref{theorem existence}}\label{sect: proof of existence}
The sets that we define next are the key of the proof, which consists of constructing laminates $\nu_j$ supported in 
\[ \bigcup_{i=1}^{j}\bigcup_{a=m_2+1}^{n-m_1-1}\left(A_{j}^{i}\cup B_{i}^{j}\cup S_{i,j}^{a}\cup S_{j,i}^{a}\right),\]
 and then approximate those laminates by homeomorphisms using Proposition \ref{prop: laminate-bisobolev}.

Given $A\in \R^{n\times n}$, let $  \sigma_1\leq\cdots\leq\sigma_n $ be its singular values. For $i,k\in\N\setminus\{0\}$ we define the following sets in the case $m_1+m_2\leq n-1$:
\begin{equation*}
A_{k}^{i}=\left\{
\begin{split}
A\in \Gamma_{+}: &|\sigma_j-(k+1)^{-1}|<\frac{(k+1)^{-2}}{4} \text{ for }j\in\{1,\ldots,n-m_1\}\text{ and }\\
&i-\frac{1}{4}< \sigma_j< i+\frac{5}{4}\text{ for }j\in\{n-m_1+1,\ldots,n\}
\end{split}
\right\},
\end{equation*}
\begin{equation*}
B_{k}^{i}=\left\{
\begin{split}
A\in \Gamma_{+}:& (k+1)^{-1}-\frac{(k+1)^{-2}}{4}< \sigma_j< k^{-1}+\frac{k^{-2}}{4}\text{ for }j\in\{1,\ldots,m_2\}\\
&\text{ and }|\sigma_j-i-1|<\frac{1}{4}\text{ for }j\in\{m_2+1,\ldots,n\}
\end{split}
\right\},
\end{equation*}
and for $a\in\{m_2+1,\ldots,n-m_1-1\}$ we define
\begin{equation*}
S_{k,i}^{a}=\left\{
\begin{split}
A\in \Gamma_{+}: &|\sigma_j-(k+1)^{-1}|<\frac{(k+1)^{-2}}{4}\text{ for }j\in\{1,\ldots,a\}\text{ and }\\
&|\sigma_j-i-1|<\frac{1}{4}\text{ for }j\in\{a+1,\ldots,n\}
\end{split}
\right\}.
\end{equation*}
We will only prove the theorem in the previous case, since, in the case $m_1+m_2\geq n$, the proof is the same using the next sets instead of the above:
\begin{equation*}
A_{k}^{i}=\left\{
\begin{split}
&A\in \Gamma_{+}: |\sigma_j-(k+1)^{-1}|<\frac{(k+1)^{-2}}{4} \text{ for }j\in\{1,\ldots,n-m_1\},\ \frac{1}{2}<\sigma_j<2\\
& \text{ for }j\in\{n-m_1+1,\ldots,m_2+1\} \text{ and } i-\frac{1}{4}< \sigma_j< i+\frac{5}{4}\text{ for }j\in\{m_2+2,\ldots,n\}
\end{split}
\right\},
\end{equation*}
and
\begin{equation*}
B_{k}^{i}=\left\{
\begin{split}
&A\in \Gamma_{+}: (k+1)^{-1}-\frac{(k+1)^{-2}}{4}< \sigma_j< k^{-1}+\frac{k^{-2}}{4}\text{ for }j\in\{1,\ldots,n-m_1-1\},\\
& \frac{1}{2}<\sigma_j<2 \text{ for }j\in\{n-m_1,\ldots,m_2\}\text{ and }|\sigma_j-i-1|<\frac{1}{4}\text{ for }j\in\{m_2+1,\ldots,n\}
\end{split}
\right\}.
\end{equation*}
The most important case is when $m_1+m_2=n-1$, where we have $S_{k,i}^{a}=\emptyset$, and, therefore, the proof is much simpler. When $m_1+m_2<n-1$ the sets $S_{k,i}^{a}$ constitute an interpolation between $A_{k}^{i}$ and $B_{k}^{i}$. We recommend the reader to focus on the case $m_1+m_2=n-1$ in a first read.

In order to write all the lemmas in a form that include all the cases, we recall the definition 
\begin{equation*}
m_i'=\begin{cases}
m_i&\text{if }m_1+m_2\leq n-1,\\
m_i-(m_1+m_2-n+1)&\text{if }m_1+m_2\geq n,
\end{cases}
\end{equation*}
for $i\in\{1,2\}$, and we define
\begin{equation*}
n'=\begin{cases}
n&\text{if }m_1+m_2\leq n-1,\\
2n-m_1-m_2-1&\text{if }m_1+m_2\geq n.
\end{cases}
\end{equation*}
In the next lemma we construct a laminate supported in $A_{k+1}^{i}\cup B_{k+1}^{i}\cup \bigcup_{a=m_2+1}^{n-m_1-1}S_{k+1,i}^{a}$ whose barycenter is a given matrix in $A_k^{i}$. We provide the proof in the case $m_1+m_2\leq n-1$. If $m_1+m_2\geq n$, we fix the eigenvalues $\{\sigma_j\}_{j=n-m_1+1}^{m_2+1} $, which are the eigenvalues in $(\frac{1}{2},2)$, and we construct the same laminate as in the first case over the other eigenvalues, i.e., given $A=\diag (\sigma_1,\ldots,\sigma_n)\in A_{k}^{i}$, let $A'=\diag (\sigma_1,\ldots,\sigma_{n-m_1},\sigma_{m_2+2},\ldots,\sigma_{n})\in\R^{n'\times n'}$ and apply Lemma \ref{lemma laminado Aki} with $n', m_1'$ and $m_2'$ to get the laminate $\nu=\sum_{\ell=1}^{N}\lambda_{\ell}\delta_{M_{\ell}'}$, where
\[M_{\ell}'=\diag (s_{1}^{\ell},\ldots,s_{n'}^{\ell}).\]
 For $\ell=1,\ldots,N$, define
 \[M_{\ell}=\diag (s_{1}^{\ell},\ldots,s_{n-m_1}^{\ell},\sigma_{n-m_1+1}(A),\ldots,\sigma_{m_2+1}(A),s_{n-m_1+1}^{\ell},\ldots,s_{n'}^{\ell}).\]
So, in the case $m_1+m_2\geq n$ we would work with the laminate $\sum_{\ell=1}^{N}\lambda_{\ell}\delta_{M_{\ell}}$.

 The bounds of $S_{k,i}^{a}$ along the paper only make sense when $m_1+m_2< n$; otherwise, the $S_{k,i}^{a}$ are empty.
\begin{lem}\label{lemma laminado Aki}
%En el paper de verdad mejor hacer que este cerca de Bki y hacer el mismo laminado
Let $A\in A_{k}^{i}$. Then there exists a laminate of finite order $\nu=\sum_{\ell=1}^{N}\lambda_{\ell}\delta_{M_{\ell}}$ such that 
\begin{enumerate}[a)]
\item\label{Aki center of mass} $\overline{\nu}=A$,
\item\label{Aki supp} $\supp (\nu)\subset A_{k+1}^{i}\cup B_{k+1}^{i}\cup \bigcup_{a=m_2+1}^{n-m_1-1}S_{k+1,i}^{a}$,
\item\label{Aki bound A} $\nu(A_{k+1}^{i})\leq 1$,
\item\label{Aki bound B} $\nu(B_{k+1}^{i})\lesssim \left(k^{2}i\right)^{m_1'+m_2'-n'}$,
\item\label{Aki bound S}  $\nu(S_{k+1,i}^{a})\lesssim\left(k^{2}i\right)^{m_1'+a-n'}$ for $a\in\{m_2+1,\ldots,n-m_1-1\}$.
\item\label{Aki distance} $M_1\in A_{k+1}^{i}$, $|A-M_1|\leq k^{-2}$, $|A^{-1}-M_{1}^{-1}|\lesssim 1$ and $1-\lambda_{1}\lesssim k^{-2}i^{-1}$.
\end{enumerate} 
\end{lem}
\begin{proof}
Since we give the proof in the case $m_1+m_2\leq n-1$, we have $m_1'=m_1$ and $m_2'=m_2$. Without loss of generality we can assume that $ A$ is a diagonal matrix, hence $A=\diag(\sigma_1,\ldots,\sigma_n)$ with $|\sigma_j-(k+1)^{-1}|<\frac{(k+1)^{-2}}{4}$ for $j\in\{1,\ldots,n-m_1\}$ and $i-\frac{1}{4}< \sigma_j< i+\frac{5}{4}$ for $j\in\{n-m_1+1,\ldots,n\}$.

Let $b=\#\{j\in\{n-m_1+1,\ldots,n\}:\sigma_j> i+\frac{3}{4}\} $.

We shall construct a family $\{B_{\ell,j}\}_{\substack{\ell=0,\ldots,n-b \\ j=0,\ldots, 2^{\ell}-1}}$ in $\Gamma_{+}$ and a family  $\{\lambda_{\ell,j}\}_{\substack{\ell=0,\ldots,n-b \\ j=0,\ldots, 2^{\ell}-1}}$ in $[0,1]$ by finite induction on $\ell$.

Let $B_{0,0}=A$, $\lambda_{0,0}=1$ and for $0\leq \ell\leq   n-b-1$,  $0\leq j\leq 2^{\ell}-1$, we assume that $\{B_{\ell,j}\}_{j=0}^{2^{\ell}-1}$ and $\{\lambda_{\ell,j}\}_{j=0}^{2^{\ell}-1}$ have been defined, $B_{\ell,j}$ are diagonal, $\lambda_{\ell,j}\geq 0$,
\begin{equation}\label{Aki barycenter}
\sum_{j=0}^{2^{\ell}-1}\lambda_{\ell,j}=1,\qquad B_{0,0}=\sum_{j=0}^{2^{\ell}-1}\lambda_{\ell,j}B_{\ell,j}
\end{equation}
\begin{equation}\label{Aki Blj es laminado}
\sum_{j=0}^{2^{\ell}-1}\lambda_{\ell,j}\delta_{B_{\ell,j}}\in\mathcal{L}(\R^{n\times n}),
\end{equation}
and when we let 
\[\beta_{\ell,j}^{1}:=\#\left\{\alpha\in\{1,\ldots,n-m_1\}: |\left(B_{\ell,j}\right)_{\alpha,\alpha}-(k+1)^{-1}|<\frac{(k+1)^{-2}}{4}\right\},\]
\[\beta_{\ell,j}^{2}:=\#\left\{\alpha\in\{1,\ldots,n-m_1\}: |\left(B_{\ell,j}\right)_{\alpha,\alpha}-(k+2)^{-1}|<\frac{(k+2)^{-2}}{4}\right\},\]
\[\beta_{\ell,j}^{3}:=\#\left\{\alpha\in\{1,\ldots,n-m_1\}: |\left(B_{\ell,j}\right)_{\alpha,\alpha}-i-1|<\frac{1}{4}\right\},\]
\[\gamma_{\ell,j}^{1}:=\#\left\{\alpha\in\{n-m_1+1,\ldots,n-b\}: i-\frac{1}{4}<\left(B_{\ell,j}\right)_{\alpha,\alpha}\leq i+\frac{3}{4}\right\},\]
\[\gamma_{\ell,j}^{2}:=\#\left\{\alpha\in\{n-m_1+1,\ldots,n-b\}: |\left(B_{\ell,j}\right)_{\alpha,\alpha}-(k+2)^{-1}|<\frac{(k+2)^{-2}}{4}\right\},\]
\[\gamma_{\ell,j}^{3}:=\#\left\{\alpha\in\{n-m_1+1,\ldots,n-b\}: |\left(B_{\ell,j}\right)_{\alpha,\alpha}-i-1|<\frac{1}{4}\right\},\]
then
\begin{equation}\label{Aki small singlevalues leq n-m1}
\beta_{\ell,j}^{1}+\beta_{\ell,j}^{2}+\gamma_{\ell,j}^{2}\leq n-m_1,
\end{equation}
\begin{equation}\label{Aki big singlevalues leq n-m2}
\beta_{\ell,j}^{3}+\gamma_{\ell,j}^{1}+\gamma_{\ell,j}^{3}\leq n-m_2-b,
\end{equation}
\begin{equation}\label{Aki all singlevalues = n}
\beta_{\ell,j}^{1}+\beta_{\ell,j}^{2}+\beta_{\ell,j}^{3}+\gamma_{\ell,j}^{1}+\gamma_{\ell,j}^{2}+\gamma_{\ell,j}^{3}=n-b,
\end{equation}
\begin{equation}\label{Aki bound lambda}
\lambda_{\ell,j}\leq\left(\frac{2}{i(k+1)^{2}}\right)^{\beta_{\ell,j}^{3}}\left(\frac{2}{i}\right)^{\gamma_{\ell,j}^{2}},
\end{equation}
\begin{equation}\label{Aki Bl0}
B_{\ell,0}=\diag\left(\underbrace{(k+2)^{-1},\ldots,(k+2)^{-1}}_{\min\{\ell,n-m_1\}},\sigma_{\min\{\ell,n-m_1\}+1},\ldots,\sigma_n\right)
\end{equation}
and
\begin{equation}\label{Aki bound lambda0}
\lambda_{\ell,0}=\prod_{j=1}^{\min\{\ell,n-m_1\}}\frac{i+1-\sigma_j}{i+1-(k+2)^{-1}}.
\end{equation}
%\begin{equation}\label{Aki det}
%\frac{\det(B_{\ell,j})}{\det(A)}\leq\left(\frac{k}{k+1}\right)^{\beta_{\ell,j}^{(k+1)^{-1}}}\left(ki\right)^{\beta_{\ell,j}^{i}}\left(\frac{1}{(k+1)(i-1)}\right)^{\gamma_{\ell,j}^{(k+1)^{-1}}}\left(\frac{i}{i-1}\right)^{\gamma_{\ell,j}^{i}}.
%\end{equation}
Moreover, for those $B_{\ell,j}\notin  A_{k+1}^{i}\cup B_{k+1}^{i}$, we have
\begin{equation}\label{Aki number of new singlevalues of Blj}
\beta_{\ell,j}^{1}+\gamma_{\ell,j}^{1}\leq n-\ell-b.
\end{equation}

Since (\ref{Aki all singlevalues = n}) holds and $m_1+m_2\leq n-1$, we see that (\ref{Aki small singlevalues leq n-m1}) and (\ref{Aki big singlevalues leq n-m2}) cannot be equalities at the same time.

We observe that the sets appearing in the definitions of $\beta_{\ell,j}^{a},\gamma_{\ell,j}^{a}$ for $a=1,2,3$ are pairwise disjoint; we also see that
\[\beta_{0,0}^{2}=\gamma_{0,0}^{2}=\beta_{0,0}^{3}=\gamma_{0,0}^{3}=0,\quad \beta_{0,0}^{1}=n-m_1,\qquad \gamma_{0,0}^{1}=m_1-b.\]
Now, we start with the induction. If $\beta_{\ell,j}^{2}+\gamma_{\ell,j}^{2}=n-m_1$, then $B_{\ell,j}\in A_{k+1}^{i}$, if $\beta_{\ell,j}^{3}+\gamma_{\ell,j}^{3}=n-m_2-b$, then $B_{\ell,j}\in B_{k+1}^{i}$, and, if $a:=\beta_{\ell,j}^{2}+\gamma_{\ell,j}^{2}<n-m_1$, $\beta_{\ell,j}^{3}+\gamma_{\ell,j}^{3}<n-m_2-b$ and $\beta_{\ell,j}^{1}+\gamma_{\ell,j}^{1}=0$ then $B_{\ell,j}\in S_{k+1,i}^{a}$.

Now, for $j=0,\ldots,2^{\ell+1}-1$ we construct $B_{\ell+1,j}$ and $\lambda_{\ell+1,j}$.

If $B_{\ell,\floor{\frac{j}{2}}}\in  A_{k+1}^{i}\cup B_{k+1}^{i}\cup \bigcup_{a=n-m_1+1}^{m_2-1}S_{k+1,i}^{a} $ we define $B_{\ell+1,j}=B_{\ell,\floor{\frac{j}{2}}}$ and 
\[\lambda_{\ell+1,j}=\begin{cases}
\lambda_{\ell,\frac{j}{2}}, &\text{if } j\text{ is even},\\
0, &\text{if } j\text{ is odd}.
\end{cases}\]
So it is clear that (\ref{Aki small singlevalues leq n-m1})--(\ref{Aki bound lambda}) are satisfied.

In the case $B_{\ell,\floor{\frac{j}{2}}}\notin  A_{k+1}^{i}\cup B_{k+1}^{i}\cup \bigcup_{a=n-m_1+1}^{m_2-1}S_{k+1,i}^{a} $, we have 
\begin{equation}\label{Aki 2+2}
\beta_{\ell,\floor{\frac{j}{2}}}^{2}+\gamma_{\ell,\floor{\frac{j}{2}}}^{2}<n-m_1,
\end{equation}
\begin{equation}\label{Aki 3+3}
\beta_{\ell,\floor{\frac{j}{2}}}^{3}+\gamma_{\ell,\floor{\frac{j}{2}}}^{3}<n-m_2-b,
\end{equation}
 and 
 \begin{equation}\label{Aki 1+1}
 \beta_{\ell,\floor{\frac{j}{2}}}^{1}+\gamma_{\ell,\floor{\frac{j}{2}}}^{1}>0,
 \end{equation}
 and, we divide the construction of $B_{\ell+1,j}$ into two cases, according to whether (\ref{Aki 1ºcaso}) or (\ref{Aki 2ºcaso}) holds. We observe that if $B_{\ell,0}\in A_{k+1}^{i}\cup B_{k+1}^{i}\cup \bigcup_{a=n-m_1+1}^{m_2-1}S_{k+1,i}^{a}$ then $B_{\ell,0}\in A_{k+1}^{i}$, which happens if and only if $\ell<n-m_1$. Hence, if $\ell\geq n-m_1$, (\ref{Aki Bl0}) and (\ref{Aki bound lambda0}) are satisfied, whereas if $\ell< n-m_1$ then (\ref{Aki 1ºcaso}) is satisfied for $j=0$.
 
 If
\begin{equation}\label{Aki 1ºcaso}
\beta_{\ell,\floor{\frac{j}{2}}}^{1}>0\text{ and }\beta_{\ell,\floor{\frac{j}{2}}}^{3}+\gamma_{\ell,\floor{\frac{j}{2}}}^{1}+\gamma_{\ell,\floor{\frac{j}{2}}}^{3}<n-m_2-b,
\end{equation}
let $\alpha\in\{1,\ldots,n-m_1\}$ be the smallest number such that $|(B_{\ell,\floor{\frac{j}{2}}})_{\alpha,\alpha}-(k+1)^{-1}|<\frac{(k+1)^{-2}}{4}$. Then, we define
\[ B_{\ell+1,j}=\begin{cases}
B_{\ell,\floor{\frac{j}{2}}}+\diag\left(\underbrace{0,\ldots,0}_{\alpha-1},(k+2)^{-1}-(B_{\ell,\floor{\frac{j}{2}}})_{\alpha,\alpha},\underbrace{0,\ldots,0}_{n-\alpha}  \right), & \text{if }  j\text{ is even},\\
B_{\ell,\floor{\frac{j}{2}}}+\diag\left(\underbrace{0,\ldots,0}_{\alpha-1},i+1-(B_{\ell,\floor{\frac{j}{2}}})_{\alpha,\alpha},\underbrace{0,\ldots,0}_{n-\alpha}  \right), & \text{if }  j\text{ is odd},
\end{cases}\]
and
\[\lambda_{\ell+1,j}=\begin{cases}
\frac{i+1-(B_{\ell,\floor{\frac{j}{2}}})_{\alpha,\alpha}}{i+1-(k+2)^{-1}}\lambda_{\ell,\floor{\frac{j}{2}}}\leq \lambda_{\ell,\floor{\frac{j}{2}}}, &\text{if } j\text{ is even},\\
\frac{(B_{\ell,\floor{\frac{j}{2}}})_{\alpha,\alpha}-(k+2)^{-1}}{i+1-(k+2)^{-1}}\lambda_{\ell,\floor{\frac{j}{2}}}\leq \frac{2}{i(k+1)^{2}}\lambda_{\ell,\floor{\frac{j}{2}}}, &\text{if } j\text{ is odd}.
\end{cases}\]
Then
\[B_{\ell,\frac{j}{2}}=\frac{i+1-(B_{\ell,\frac{j}{2}})_{\alpha,\alpha}}{i+1-(k+2)^{-1}} B_{\ell+1,j}+\frac{(B_{\ell,\frac{j}{2}})_{\alpha,\alpha}-(k+2)^{-1}}{i+1-(k+2)^{-1}}B_{\ell+1,j+1}\text{ for }j\in\{0,\ldots,2^{\ell+1}-1\}\text{ even},\]
and hence
\[\sum_{j=0}^{2^{\ell+1}-1}\lambda_{\ell+1,j}\delta_{B_{\ell+1,j}}\in\mathcal{L}(\R^{n\times n}),\]
\[\beta_{\ell+1,j}^{1}=\beta_{\ell,\floor{\frac{j}{2}}}^{1}-1,\]
\[\beta_{\ell+1,j}^{2}=\begin{cases}
\beta_{\ell,\floor{\frac{j}{2}}}^{2}+1,&\text{if } j\text{ is even},\\
\beta_{\ell,\floor{\frac{j}{2}}}^{2}, &\text{if } j\text{ is odd},
\end{cases}\]
\[\beta_{\ell+1,j}^{3}=\begin{cases}
\beta_{\ell,\floor{\frac{j}{2}}}^{3},&\text{if } j\text{ is even},\\
\beta_{\ell,\floor{\frac{j}{2}}}^{3}+1, &\text{if } j\text{ is odd},
\end{cases}\]
\[\gamma_{\ell+1,j}^{1}=\gamma_{\ell,\floor{\frac{j}{2}}}^{1},\quad \gamma_{\ell+1,j}^{2}=\gamma_{\ell,\floor{\frac{j}{2}}}^{2}\quad\text{ and }\gamma_{\ell+1,j}^{3}=\gamma_{\ell,\floor{\frac{j}{2}}}^{3}.\]
Therefore, (\ref{Aki barycenter})--(\ref{Aki number of new singlevalues of Blj}) are satisfied for $\ell+1$.

If
\begin{equation}\label{Aki 2ºcaso}
\beta_{\ell,\floor{\frac{j}{2}}}^{1}=0\text{ or }\beta_{\ell,\floor{\frac{j}{2}}}^{3}+\gamma_{\ell,\floor{\frac{j}{2}}}^{1}+\gamma_{\ell,\floor{\frac{j}{2}}}^{3}= n-m_2-b
\end{equation}
instead of (\ref{Aki 1ºcaso}), we claim that then we have $\gamma_{\ell,\floor{\frac{j}{2}}}^{1}>0$ and 
\[\beta_{\ell,\floor{\frac{j}{2}}}^{1}+\beta_{\ell,\floor{\frac{j}{2}}}^{2}+\gamma_{\ell,\floor{\frac{j}{2}}}^{2}< n-m_1.\]
Indeed, if $\beta_{\ell,\floor{\frac{j}{2}}}^{1}=0$, it is clear thanks to (\ref{Aki 2+2}) and (\ref{Aki 1+1}), and if 
\[\beta_{\ell,\floor{\frac{j}{2}}}^{3}+\gamma_{\ell,\floor{\frac{j}{2}}}^{1}+\gamma_{\ell,\floor{\frac{j}{2}}}^{3}= n-m_2-b,\] 
we have 
\[\beta_{\ell,\floor{\frac{j}{2}}}^{1}+\beta_{\ell,\floor{\frac{j}{2}}}^{2}+\gamma_{\ell,\floor{\frac{j}{2}}}^{2}= m_2<n-m_1,\] 
and by (\ref{Aki 3+3}) we obtain $\gamma_{\ell,\floor{\frac{j}{2}}}^{1}>0$. Therefore there exists $\alpha\in\{n-m_1+1,\ldots,n-b\}$ such that 
 \[i-\frac{1}{4}<(B_{\ell,\floor{\frac{j}{2}}})_{\alpha,\alpha}\leq i+\frac{3}{4}.\]
  Then, we define
\[ B_{\ell+1,j}=\begin{cases}
B_{\ell,\floor{\frac{j}{2}}}+\diag\left(\underbrace{0,\ldots,0}_{\alpha-1},(k+2)^{-1}-(B_{\ell,\floor{\frac{j}{2}}})_{\alpha,\alpha},\underbrace{0,\ldots,0}_{n-\alpha}  \right), & \text{if }  j\text{ is even},\\
B_{\ell,\floor{\frac{j}{2}}}+\diag\left(\underbrace{0,\ldots,0}_{\alpha-1},i+1-(B_{\ell,\floor{\frac{j}{2}}})_{\alpha,\alpha},\underbrace{0,\ldots,0}_{n-\alpha}  \right), & \text{if }  j\text{ is odd},
\end{cases}\]
\[\lambda_{\ell+1,j}=\begin{cases}
\frac{i+1-(B_{\ell,\floor{\frac{j}{2}}})_{\alpha,\alpha}}{i+1-(k+2)^{-1}}\lambda_{\ell,\floor{\frac{j}{2}}}\leq\frac{2}{i}\lambda_{\ell,\floor{\frac{j}{2}}} , &\text{if } j\text{ is even},\\
\frac{(B_{\ell,\floor{\frac{j}{2}}})_{\alpha,\alpha}-(k+2)^{-1}}{i+1-(k+2)^{-1}}\lambda_{\ell,\floor{\frac{j}{2}}}\leq \lambda_{\ell,\floor{\frac{j}{2}}}, &\text{if } j\text{ is odd}.
\end{cases}\]
Then 
\[\sum_{j=0}^{2^{\ell+1}-1}\lambda_{\ell+1,j}\delta_{B_{\ell+1,j}}\in\mathcal{L}(\R^{n\times n}),\]
\[\beta_{\ell+1,j}^{1}=\beta_{\ell,\floor{\frac{j}{2}}}^{1},\quad\beta_{\ell+1,j}^{2}=\beta_{\ell,\floor{\frac{j}{2}}}^{2}, \quad\beta_{\ell+1,j}^{3}=\beta_{\ell,\floor{\frac{j}{2}}}^{3},\]
\[\gamma_{\ell+1,j}^{1}=\gamma_{\ell,\floor{\frac{j}{2}}}^{1}-1,\]
\[\gamma_{\ell+1,j}^{2}=\begin{cases}
\gamma_{\ell,\floor{\frac{j}{2}}}^{2}+1, &\text{if } j\text{ is even},\\
\gamma_{\ell,\floor{\frac{j}{2}}}^{2}, &\text{if } j\text{ is odd,}
\end{cases}\]
and
\[\gamma_{\ell+1,j}^{3}=\begin{cases}
\gamma_{\ell,\floor{\frac{j}{2}}}^{3}, &\text{if } j\text{ is even},\\
\gamma_{\ell,\floor{\frac{j}{2}}}^{3}+1, &\text{if } j\text{ is odd}.
\end{cases}\]
Therefore, (\ref{Aki barycenter}), (\ref{Aki Blj es laminado}), (\ref{Aki small singlevalues leq n-m1}), (\ref{Aki big singlevalues leq n-m2}), (\ref{Aki all singlevalues = n}), (\ref{Aki bound lambda}) and (\ref{Aki number of new singlevalues of Blj}) are satisfied for $\ell+1$.

Here ends the inductive construction of $\{B_{\ell,j}\}_{\substack{\ell=0,\ldots,n-b \\ j=0,\ldots, 2^{\ell}-1}}$. With this, we define
\[N=2^{n-b},\quad \lambda_{\ell}=\lambda_{n-b,\ell-1},\quad M_{\ell}=B_{n-b,\ell-1}\text{ for }\ell=1,\ldots,2^{n-b}\quad\text{and }\nu:=\sum_{j=1}^{2^{n-b}}\lambda_{j}\delta_{M_{j}},\]
which is a laminate by (\ref{Aki Blj es laminado}), and by (\ref{Aki all singlevalues = n}) and (\ref{Aki number of new singlevalues of Blj}) it is supported in $A_{k+1}^{i}\cup B_{k+1}^{i}\cup \bigcup_{a=m_2+1}^{n-m_1-1}S_{k+1,i}^{a}$. We shall check properties \emph{\ref{Aki center of mass})}--\emph{\ref{Aki distance})}. Property \emph{\ref{Aki center of mass})} comes from (\ref{Aki barycenter}), \emph{\ref{Aki supp})} is by (\ref{Aki number of new singlevalues of Blj}), and \emph{\ref{Aki bound A})} is obvious.

Now, we use (\ref{Aki bound lambda}) to bound the mass of $\nu$ in the different sets.

Since $\gamma_{n-b,j}^{3}\leq m_1-b$ and the matrices $B_{n-b,j}$ in $B_{k+1}^{i}$ are those such that $\beta_{n-b,j}^{3}+\gamma_{n-b,j}^{3}=n-m_2-b$, we have $\beta_{n-b,j}^{3}\geq n-m_1-m_2$ and
\[\nu(B_{k+1}^{i})=\sum_{j: \beta_{n-b,j}^{3}+\gamma_{n-b,j}^{3}=n-m_2-b }\lambda_{n-b,j}\lesssim \left(\frac{2}{i(k+1)^2}\right)^{n-m_1-m_2}\lesssim \left(k^{2}i\right)^{m_1+m_2-n}.\]
So \emph{\ref{Aki bound B})} is proved. Now we use that $\sum_{l=1}^{3}\beta_{n-b,j}^{l}=n-m_1$ and that for $a\in\{m_2+1,\ldots,n-m_1-1\}$ the matrices $B_{n-b,j}$ in $S_{k+1,i}^{a}$ are those such that $\beta_{n-b,j}^{2}+\gamma_{n-b,j}^{2}=a$, $\beta_{n-b,j}^{1}=\gamma_{n-b,j}^{1}=0$, so
\[\beta_{n-b,j}^{3}= n-m_1- \beta_{n-b,j}^{2}=n-m_1+ \gamma_{n-b,j}^{2}-a\geq n-m_1-a,\]
 and, hence,
\[\nu(S_{k+1,i}^{a})=\sum_{j: \beta_{n-b,j}^{2}+\gamma_{n-b,j}^{2}=a }\lambda_{n-b,j}\lesssim \left(\frac{2}{i(k+1)^2}\right)^{n-m_1-a}\lesssim \left(k^{2}i\right)^{m_1+a-n}.\]
Therefore we have \emph{\ref{Aki bound S})}. Finally, we prove \emph{\ref{Aki distance})}. Thanks to (\ref{Aki Bl0}) we get
\[M_1=\diag \left(\underbrace{(k+2)^{-1},\ldots,(k+2)^{-1}}_{n-m_1},\sigma_{n-m_1+1},\ldots,\sigma_n\right)\]
and due to (\ref{Aki bound lambda0}) we obtain
\[\lambda_1=\prod_{j=1}^{n-m_1}\frac{i+1-\sigma_j}{i+1-(k+2)^{-1}}.\]
Therefore, using $k^{-1}-(k+1)^{-1}\leq k^{-2}$ we get
\[|A-M_1|\leq k^{-2},\quad |A^{-1}-M^{-1}|\lesssim 1\]
and, noting that $n-m_1\geq 2$ we also obtain
\begin{align*}
1-\lambda_1&\leq 1-\left(\frac{i+1-(k+1)^{-1}-\frac{(k+1)^{-2}}{4}}{i+1-(k+2)^{-1}}\right)^{n-m_1}\\
&\lesssim \frac{\left((k+1)^{-1}+\frac{(k+1)^{-2}}{4}\right)^{n-m_1-1}-(k+2)^{m_1-n+1}}{i}\lesssim k^{-2}i^{-1}.
\end{align*}
Hence, the proof is finished.
\end{proof}
In the next two lemmas we give the laminates starting in $B_k^{i}$ and $S_{k,i}^{a}$. We will not show them, since their construction mimic that of Lemma \ref{lemma laminado Aki}. The main difference in the proofs of Lemmas \ref{lemma laminado Bki} and \ref{lemma laminado Ski} with the one of Lemma \ref{lemma laminado Aki} is that we split the eigenvalues $\sigma$ close to $i+1$ in $(k+1)^{-1}$ and $\sigma+1$. Also, in the proof of Lemma \ref{lemma laminado Bki} we start splitting the eigenvalues close to $i+1$ and we do not split the eigenvalues in $\left((k+1)^{-1}-\frac{(k+1)^{-2}}{4},(k+1)^{-1}+\frac{(k+1)^{-2}}{4}\right)$.
\begin{lem}\label{lemma laminado Bki}
Let $A\in B_{k}^{i}$. Then there exists a laminate of finite order $\nu=\sum_{\ell=1}^{N}\lambda_{\ell}\delta_{M_{\ell}}$ such that 
\begin{itemize}
\item $\overline{\nu}=A$,
\item $\supp (\nu)\subset A_{k}^{i+1}\cup B_{k}^{i+1}\cup \bigcup_{a=m_2+1}^{n-m_1-1}S_{k,i+1}^{a}$,
\item $\nu(A_{k}^{i+1})\lesssim i^{m_1'+m_2'-n'}$,
\item $\nu(B_{k}^{i+1})- \left(\frac{i+2}{i+3}\right)^{n'-m_2'}\lesssim (ki)^{-2}$,
\item  $\nu(S_{k,i+1}^{a})\lesssim i^{m_2'-a}$ for $a\in\{m_2+1,\ldots,n-m_1-1\}$,
\item\label{Bki distance} $M_1\in B_{k}^{i+1}$, $|A-M_1|=1$, $|A^{-1}-M_{1}^{-1}|\leq i^{-2}$ and $1-\lambda_{1}\lesssim i^{-1}$.
\end{itemize} 
\end{lem}
\begin{lem}\label{lemma laminado Ski}
Let $a_0\in\{m_2+1,\ldots,n-m_1-1\}$ and $A\in S_{k,i}^{a_0}$. Then there exists a laminate of finite order $\nu=\sum_{\ell=1}^{N}\lambda_{\ell}\delta_{M_{\ell}}$ such that 
\begin{itemize}
\item $\overline{\nu}=A$,
\item $\supp (\nu)\subset A_{k+1}^{i+1}\cup B_{k+1}^{i+1}\cup \bigcup_{a=m_2+1}^{n-m_1-1}S_{k+1,i+1}^{a}$,
\item $\nu(A_{k+1}^{i+1})\lesssim i^{a_0+m_1-n}$,
\item $\nu(B_{k+1}^{i+1})\lesssim (k^{2}i)^{m_2-a_0}$,
\item $\nu(S_{k+1,i+1}^{a})\lesssim (k^{2}i)^{a-a_0}$ if $a\in\{m_2+1,\ldots,a_0-1\}$,
\item $ \nu(S_{k+1,i+1}^{a_0})- \left(\frac{i+2}{i+3}\right)^{n-a_0}\lesssim (ki)^{-2}$,
\item $\nu(S_{k+1,i+1}^{a})\lesssim i^{a_0-a}$ if $a\in\{a_0+1,\ldots,n-m_1-1\}$,
\item  $\sum_{\ell=1}^{N}\lambda_{\ell}|A-M_{\ell}|\lesssim 1$,
\item $\frac{1}{\det (A)}\sum_{\ell=1}^{N}\lambda_{\ell}\det (M_{\ell})|A^{-1}-M_{\ell}^{-1}|\lesssim 1$.
\end{itemize} 
\end{lem}
In the proof of the last lemma, besides adapting the proof of Lemma \ref{lemma laminado Aki}, we follow the proof of \emph{\ref{Aji conv f})} and \emph{\ref{Aji conv f-1})} of Lemma \ref{lemma laminado Aji} below to show the last two items.

In the next lemma we put together Lemmas \ref{lemma laminado Aki}, \ref{lemma laminado Bki} and \ref{lemma laminado Ski} to construct laminates whose support is in the set in which we are interested. Again, all the bounds of $S_{j+1,i}^{a}$ have sense if and only if $m_1+m_2\leq n-1$; otherwise these sets are empty. 
\begin{lem}\label{lemma laminado Aji}
Let $i,j\in\N$, $i\leq j$, and $A\in A_j^{i}$. Then there exists a laminate $\nu=\sum_{k=1}^{N}\lambda_{k}\delta_{M_{k}}$ such that 
\begin{enumerate}[a)]
\item $\overline{\nu}=A$,
\item $\supp (\nu)\subset \bigcup_{b=0}^{j-i+1}A_{j+1}^{i+b}\cup B_{j+1}^{j+1}\cup \bigcup_{b=0}^{j-i+1}\bigcup_{a=m_2+1}^{n-m_1-1}S_{j+1,i+b}^{a}$,
\item $\nu(A_{j+1}^{i})\leq 1 $,
\item $\nu(A_{j+1}^{i+b})\lesssim j^{2(m_1'+m_2'-n')}  \frac{i^{m_1'}}{(i+b)^{m_1'+2} }$, for $b\in\{1,\ldots,j-i+1\}$,
\item $\nu(B_{j+1}^{j+1})\lesssim j^{2m_1'+3m_2'-3n'}i^{m_1'}$,
\item  $\nu(S_{j+1,i}^{a})\lesssim (j^{2} i)^{m_1'+a-n'}$,
\item  $\nu(S_{j+1,i+b}^{a})\lesssim j^{2(m_1'+m_2'-n')}i^{m_1'}(i+b)^{2m_2'-a-n'}$ for $a\in\{m_2+1,\ldots,n-m_1-1\}$ and $b\in\{1,\ldots,j-i+1\}$,
\item\label{Aji conv f}  $\sum_{k=1}^{N}\lambda_{k}|A-M_{k}|\lesssim j^{-2}$,
\item\label{Aji conv f-1} $\frac{1}{\det (A)}\sum_{k=1}^{N}\lambda_{k}\det (M_{k})|A^{-1}-M_{k}^{-1}|\lesssim 1$.
\end{enumerate}
\end{lem}
\begin{proof}
Let $C$ be a constant bigger than those in Lemma \ref{lemma laminado Bki}. It is enough construct a sequence of laminates $\{\nu_{\ell}\}_{\ell=1}^{j-i+2}$ such that
\begin{enumerate}[1)]
\item\label{pro Aji barycenter} $\overline{\nu}_{\ell}=A$,
\item $\supp (\nu_{\ell})\subset \bigcup_{b=0}^{\ell-1}A_{j+1}^{i+b}\cup B_{j+1}^{i+\ell-1}\cup \bigcup_{b=0}^{\ell-1}\bigcup_{a=m_2+1}^{n-m_1-1}S_{j+1,i+b}^{a}$,
\item $\nu_{\ell}(A_{j+1}^{i+b})\leq j^{2(m_1'+m_2'-n')}  \frac{i^{m_1'}}{(i+b)^{m_1'+2} } $, for $b\in\{1,\ldots,\ell-1\}$,
\item $\nu_{\ell}(B_{j+1}^{i+\ell-1})\lesssim j^{2(m_1'+m_2'-n')}i^{m_1'}(i+\ell+1)^{m_2'-n'}\left( 1+j^{-2}\sum_{k=1}^{\ell-1}2^{n}C(i+k)^{-2} \right)$,
\item  $\nu_{\ell}(S_{j+1,i}^{a})\lesssim (j^{2} i)^{m_1'+a-n'}$,
\item  $\nu_{\ell}(S_{j+1,i+b}^{a})\lesssim j^{2(m_1'+m_2'-n')}i^{m_1'}(i+b)^{2m_2'-a-n'}$ for $a\in\{m_2+1,\ldots,n-m_1-1\}$ and $b\in\{1,\ldots,\ell-1\}$,
\item\label{pro Aji exist M} $M_{1}\in A_{j+1}^{i}\cap \supp(\nu_{\ell})$ such that $|A-M_{1}|\leq j^{-2}$, $|A^{-1}-M_{1}^{-1}|\lesssim 1$, $1-\nu_{\ell}(M_{1})\lesssim j^{-2}i^{-1}$, with $M_1$ being the one of Lemma \ref{lemma laminado Aki},
\end{enumerate}
and prove later \emph{\ref{Aji conv f})} and \emph{\ref{Aji conv f-1})}.

Let $\nu_{1}$ be the laminate of Lemma \ref{lemma laminado Aki}, which satisfies
\begin{itemize}
\item $\overline{\nu}_{1}=A$,
\item $\supp (\nu_{1})\subset A_{j+1}^{i}\cup B_{j+1}^{i}\cup \bigcup_{a=m_2+1}^{n-m_1-1}S_{j+1,i}^{a}$,
\item $\nu_{1}(B_{j+1}^{i})\lesssim \left(j^{2}i\right)^{m_1'+m_2'-n'}$,
\item  $\nu_{1}(S_{j+1,i}^{a})\lesssim\left(j^{2}i\right)^{m_1'+a-n'}$ for $a\in\{m_2+1,\ldots,n-m_1-1\}$,
\item $\exists M\in A_{j+1}^{i}\cap \supp(\nu_{1})$ such that $|A-M|\leq j^{-2}$, $|A^{-1}-M^{-1}|\lesssim 1$ and $1-\nu_{1}(M)\lesssim j^{-2}i^{-1}$.
\end{itemize}
Therefore, denoting by $M_1$ the matrix $M$, we get that $\nu_{1}$ satisfies \ref{pro Aji barycenter})--\ref{pro Aji exist M}).

Now, we proceed by induction and assume that $\nu_{\ell}$ has been constructed with the properties \ref{pro Aji barycenter})--\ref{pro Aji exist M}). For each $B\in \supp(\nu_{\ell})\cap B_{j+1}^{i+\ell-1}$, let $\nu_B$ be the laminate given by Lemma \ref{lemma laminado Bki}, which satisfies
\begin{itemize}
\item $\overline{\nu}_B=B$,
\item $\supp (\nu_B)\subset A_{j+1}^{i+\ell}\cup B_{j+1}^{i+\ell}\cup \bigcup_{a=m_2+1}^{n-m_1-1}S_{j+1,i+\ell}^{a}$,
\item $\nu_B(A_{j+1}^{i+\ell})\lesssim (i+\ell)^{m_1'+m_2'-n'}$,
\item $\nu_B(B_{j+1}^{i+\ell})- \left(\frac{i+\ell+1}{i+\ell+2}\right)^{n'-m_2'}\leq C (i+\ell)^{-2}j^{-2}$,
\item  $\nu_B(S_{j+1,i+\ell}^{a})\lesssim (i+\ell)^{m_2'-a}$ for $a\in\{m_2+1,\ldots,n-m_1-1\}$.
\end{itemize} 
We define 
\[\nu_{\ell+1}=\nu_{\ell}+\sum_{B\in \supp(\nu_{\ell})\cap B_{j+1}^{i+\ell-1}}\nu_{\ell}(B)(\nu_B-\delta_B).\]
Thanks to Corollary \ref{cor: laminate}, $\nu_{\ell+1}$ is a laminate. Moreover, it is clear that $\overline{\nu}_{\ell+1}=A$ and
\[\supp (\nu_{\ell+1})\subset \bigcup_{b=0}^{\ell}A_{j+1}^{i+b}\cup B_{j+1}^{i+\ell}\cup \bigcup_{b=0}^{\ell}\bigcup_{a=m_2+1}^{n-m_1-1}S_{j+1,i+b}^{a}.\]

Now, observe that the matrices in $(A_{j+1}^{i+\ell}\cup B_{j+1}^{i+\ell}\cup\bigcup_{a=m_2+1}^{n-m_1-1}S_{j+1,i+\ell}^{a})\cap\supp (\nu_{\ell+1})$ are those in the support of $\nu_B$ for $B$ in $B_{j+1}^{i+\ell-1}\cap\supp (\nu_{\ell})$. Therefore, using  $m_1'+m_2'\leq n'-1$ we have
\begin{align*}
&\nu_{\ell+1}(A_{j+1}^{i+\ell})= \sum_{B\in B_{j+1}^{i+\ell-1}}\nu_{\ell}(B)\nu_B(A_{j+1}^{i+\ell})\leq \sum_{B\in B_{j+1}^{i+\ell-1}}\nu_{\ell}(B)(i+\ell)^{m_1'+m_2'-n'}\\
&= \nu_{\ell}(B_{j+1}^{i+\ell-1})(i+\ell)^{m_1'+m_2'-n'}\\
&\lesssim j^{2(m_1'+m_2'-n')}i^{m_1'}(i+\ell+1)^{m_2'-n'}\left( 1+j^{-2}\sum_{k=1}^{\ell-1}2^{n}C(i+k)^{-2} \right)(i+\ell)^{m_1'+m_2'-n'}\\
&\leq j^{2(m_1'+m_2'-n')}i^{m_1'}(i+\ell+1)^{-m_1'-2}\left( 1+j^{-2}\sum_{k=1}^{\ell-1}2^{n}C(i+k)^{-2} \right),
\end{align*}
Therefore
\[\nu_{\ell+1}(A_{j+1}^{i+\ell})\lesssim j^{2(m_1'+m_2'-n')} \frac{i^{m_1'}}{(i+\ell)^{m_1'+2}}. \]
For the bound of $\nu_{\ell+1}(B_{j+1}^{i+\ell})$ we need the following estimate, in which we use $j\geq 2^{n}C$:
\begin{align*}
&\left( 1+j^{-2}\sum_{k=1}^{\ell-1}2^{n}C(i+k)^{-2} \right)\left(1+C\left(\frac{i+\ell+2}{i+\ell+1}\right)^{n'-m_2'}j^{-2}(i+\ell)^{-2} \right)\\
&\leq 1+j^{-2}\sum_{k=1}^{\ell-1}2^{n}C(i+k)^{-2}+2C2^{n-1}j^{-2}(i+\ell)^{-2} \leq 1+j^{-2}\sum_{k=1}^{\ell}2^{n}C(i+k)^{-2}.
\end{align*}
Proceeding in the same way as in the bound of $\nu_{\ell+1}(A_{j+1}^{i+\ell})$, we obtain the following bound. Note that the constant corresponding to the symbol $\lesssim$ is the same for all $\ell=1,\ldots, j-i+2$:
\begin{align*}
&\nu_{\ell+1}(B_{j+1}^{i+\ell})=\sum_{B\in B_{j+1}^{i+\ell-1}}\nu_{\ell}(B)\nu_B(B_{j+1}^{i+\ell})\leq \sum_{B\in B_{j+1}^{i+\ell-1}}\nu_{\ell}(B)\left[ \left(\frac{i+\ell+1}{i+\ell+2}\right)^{n'-m_2'}+C(i+\ell)^{-2}j^{-2} \right]\\
&= \nu_{\ell}(B_{j+1}^{i+\ell-1})\left[ \left(\frac{i+\ell+1}{i+\ell+2}\right)^{n'-m_2'}+C(i+\ell)^{-2}j^{-2} \right]\\
&\lesssim j^{2(m_1'+m_2'-n')}i^{m_1'}(i+\ell+1)^{m_2'-n'}\left( 1+j^{-2}\sum_{k=1}^{\ell-1}2^{n}C(i+k)^{-2} \right)\\
&\quad\times\left[ \left(\frac{i+\ell+1}{i+\ell+2}\right)^{n'-m_2'}+C(i+\ell)^{-2}j^{-2} \right]\\
&\leq j^{2(m_1'+m_2'-n')}i^{m_1'}(i+\ell+2)^{m_2'-n'}\left( 1+j^{-2}\sum_{k=1}^{\ell}2^{n}C(i+k)^{-2} \right).
\end{align*}
Next we bound $\nu_{\ell+1}(S_{j+1,i+\ell}^{a})$ for $a\in\{m_2+1,\ldots,n-m_1-1\}$:
\begin{align*}
&\nu_{\ell+1}(S_{j+1,i+\ell}^{a})= \sum_{B\in B_{j+1}^{i+\ell-1}}\nu_{\ell}(B)\nu_B(S_{j+1,i+\ell}^{a})\leq \sum_{B\in B_{j+1}^{i+\ell-1}}\nu_{\ell}(B)(i+\ell)^{m_2'-a}\\
&= \nu_{\ell}(B_{j+1}^{i+\ell-1})(i+\ell)^{m_2'-a}\\
&\lesssim  j^{2(m_1'+m_2'-n')}i^{m_1'}(i+\ell+1)^{m_2'-n'}\left( 1+j^{-2}\sum_{k=1}^{\ell-1}2^{n}C(i+k)^{-2} \right)(i+\ell)^{m_2'-a}\\
&\lesssim j^{2(m_1'+m_2'-n')}i^{m_1'}(i+\ell)^{2m_2'-a-n'}.
\end{align*}
For $b\in\{0,\ldots,\ell\}$, we also have that 
\[\nu_{\ell+1}\lfloor_{A_{j+1}^{i+b}}=\nu_{\ell}\lfloor_{A_{j+1}^{i+b}}\]
and
\[\nu_{\ell+1}\lfloor_{S_{j+1,i+b}^{a}}=\nu_{\ell}\lfloor_{S_{j+1,i+b}^{a}}\text{ for }a\in\{m_2+1,\ldots,n-m_1-1\}.\]
Here, $\lfloor$ denotes the restriction of a measure. Therefore, $\nu_{\ell+1}$ satisfies \ref{pro Aji barycenter})--\ref{pro Aji exist M}). Here ends the inductive construction of $\{\nu_{\ell}\}_{\ell=1}^{j-i+2}$.

Now, we define $\nu=\nu_{j-i+2}=\sum_{k=1}^{N}\lambda_{k}\delta_{M_{\ell}}$, so $\lambda_1=\nu_{1}(M_1)$ and we recall that
\[M_1\in A_{j+1}^{i}\cap \supp(\nu),\, |A-M_1|\leq j^{-2},\, |A^{-1}-M_{1}^{-1}|\lesssim 1\text{ and }1-\lambda_1\lesssim j^{-2}i^{-1}.\]
Finally we have to prove \emph{\ref{Aji conv f})} and \emph{\ref{Aji conv f-1})}.  To show \emph{\ref{Aji conv f})} we need the following estimate of the distance between the matrices in the support of $\nu$ and $A$:
\begin{equation*}
|A-M|\leq |A|+|M|\lesssim\begin{cases}
i+b&\text{if }M\in A_{j+1}^{i+b}\setminus \{M_1\}\text{ for some }b\in\{0,\ldots,j-i+1\},\\
j&\text{if }M\in B_{j+1}^{j+1},\\
i+b&\text{if }M\in S_{j+1,i+b}^{a}\text{ for some }b\in\{0,\ldots,j-i+1\}.
\end{cases}
\end{equation*}
Now we split the sum $\sum_{k=1}^{N}\lambda_{k}|A-M_{k}|$ over the different sets and we bound those sums:
\[\lambda_1|A-M_1|\lesssim j^{-2},\]
\begin{align*}
\sum_{k: M_k\in A_{j+1}^{i}\setminus \{M_1\}}\lambda_{k}|A-M_{k}|\lesssim\sum_{k: M_k\in A_{j+1}^{i}\setminus \{M_1\}}\lambda_{k} i\leq i(1-\lambda_1)\lesssim j^{-2}.
\end{align*}
  In the following estimate we use $\sum_{b=0}^{\infty}(i+b)^{-m_1-1}\lesssim i^{-m_1}$ and $m_1'+m_2'-n'\leq -1$:
  \begin{align*}
&\sum_{b=1}^{j-i+1}\sum_{k: M_k\in A_{j+1}^{i+b}}\lambda_{k}|A-M_{k}|\lesssim\sum_{b=1}^{j-i+1}\sum_{k: M_k\in A_{j+1}^{i+b}}\lambda_{k} (i+b)= \sum_{b=1}^{j-i+1}\nu(A_{j+1}^{i+b})(i+b)\\
&\lesssim \sum_{b=1}^{j-i+1}j^{2(m_1'+m_2'-n')}  \frac{i^{m_1'}}{(i+b)^{m_1'+1} }\lesssim j^{-2}.
\end{align*}
Using that $i\leq j$ and $m_1'+m_2'\leq n'-1$, we have
\begin{align*}
&\sum_{k: M_k\in B_{j+1}^{j+1}}\lambda_{k}|A-M_{k}|\lesssim \sum_{k: M_k\in B_{j+1}^{j+1}}\lambda_{k}j=\nu(B_{j+1}^{j+1})j\lesssim j^{2m_1'+3m_2'-3n'+1}i^{m_1'}\leq j^{-2}.
\end{align*}
 In the same way as before, we bound the sum over the sets $S_{j+1,i}^{a}$ and $S_{j+1,i}^{a+b}$; using that in this case $m_1'=m_1$ and $n'=n$,  we get
  \begin{align*}
&\sum_{a=m_2+1}^{n-m_1-1}\sum_{k: M_k\in S_{j+1,i}^{a}}\lambda_{k}|A-M_{k}|\lesssim\sum_{a=m_2+1}^{n-m_1-1}\sum_{k: M_k\in S_{j+1,i}^{a}}\lambda_{k} i= \sum_{a=m_2+1}^{n-m_1-1}\nu(S_{j+1,i}^{a})i\\
&\lesssim \sum_{a=m_2+1}^{n-m_1-1}j^{2(m_1'+a-n')} i ^{m_1'+a-n'+1}\lesssim j^{-2}
\end{align*}
 and
  \begin{align*}
&\sum_{b=1}^{j-i+1}\sum_{a=m_2+1}^{n-m_1-1}\sum_{k: M_k\in S_{j+1,i+b}^{a}}\lambda_{k}|A-M_{k}|\lesssim\sum_{b=1}^{j-i+1}\sum_{a=m_2+1}^{n-m_1-1}\sum_{k: M_k\in S_{j+1,i+b}^{a}}\lambda_{k} (i+b)\\
&= \sum_{b=1}^{j-i+1}\sum_{a=m_2+1}^{n-m_1-1}\nu(S_{j+1,i+b}^{a})(i+b)\lesssim \sum_{b=1}^{j-i+1}\sum_{a=m_2+1}^{n-m_1-1}j^{2(m_1'+m_2'-n')} i ^{m_1'}(i+b)^{2m_2'-a-n'+1}\\
&\lesssim j^{-2}i^{m_1'}  \sum_{b=1}^{j-i+1}(i+b)^{m_2'-n'}\leq j^{-2}i^{m_1'}  \sum_{b=1}^{j-i+1}(i+b)^{-m_1'-1}\lesssim j^{-2}.
\end{align*}
Therefore, putting together all the previous bounds we obtain 
\[\sum_{k=1}^{N}\lambda_{k}|A-M_{k}|\lesssim j^{-2},\]
and, hence, \emph{\ref{Aji conv f})} is proved.

Now, to prove \emph{\ref{Aji conv f-1})} we need to bound the distance between the inverses. For $M\in\supp(\nu)\setminus \{M_1\}$ we have
\begin{equation*}
|A^{-1}-M^{-1}|\leq |A^{-1}|+|M^{-1}|\lesssim j.
\end{equation*}
We also need the following bound of the determinants. Since $A\in A_{j}^{i}$ and $m_1-n=m_1'-n'$ we have
\[\det(A)\gtrsim j^{m_1-n}i^{m_1'}=j^{m_1'-n'}i^{m_1'}.\]
Looking at the definition of the sets we also get
\begin{equation*}
\det(M)\lesssim
\begin{cases}
j^{m_1'-n'}(i+b)^{m_1'}&\text{if }M\in A_{j+1}^{i+b}\text{ for some }b\in\{0,\ldots,j-i+1\},\\
j^{n'-2m_2'}&\text{if }M\in B_{j+1}^{j+1},\\
j^{-a}(i+b)^{n'-a}&\text{if }M\in S_{j+1,i+b}^{a}\text{ for some }b\in\{0,\ldots,j-i+1\}.
\end{cases}
\end{equation*}
Hence we get
\begin{equation*}
\frac{\det(M)}{\det (A)}|A^{-1}-M^{-1}|\lesssim
\begin{cases}
j(i+b)^{m_1'}i^{-m_1'}&\text{if }M\in A_{j+1}^{i+b}\setminus \{M_1\}\text{ for some }\\
&b\in\{0,\ldots,j-i+1\},\\
j^{2n'-2m_2'-m_1'+1}i^{-m_1'}&\text{if }M\in B_{j+1}^{j+1},\\
j^{n'-m_1'-a+1}(i+b)^{n'-a}i^{-m_1'}&\text{if }M\in S_{j+1,i+b}^{a}\text{ for some }\\
&b\in\{0,\ldots,j-i+1\}.
\end{cases}
\end{equation*}
Now we split the sum $\sum_{k=1}^{N}\lambda_{k}\frac{\det(M_k)}{\det(A)}|A^{-1}-M_{k}^{-1}|$ over the different sets and bound those sums:
\[\lambda_1\frac{\det(M_1)}{\det(A)}|A^{-1}-M_{1}^{-1}|\lesssim \frac{\det(M_1)}{\det(A)}\lesssim j^{m_1'-n'}i^{m_1'}j^{n'-m_1'}i^{-m_1'}=1,\]
\begin{align*}
\sum_{k: M_k\in A_{j+1}^{i}\setminus \{M_1\}}\lambda_{k}\frac{\det(M_k)}{\det(A)}|A^{-1}-M_{k}^{-1}|\lesssim\sum_{k: M_k\in A_{j+1}^{i}\setminus \{M_1\}}\lambda_{k} j\leq j(1-\lambda_1)\lesssim j^{-1}i^{-1}\leq 1.
\end{align*}
  In the following estimate we use $\sum_{b=0}^{\infty}(i+b)^{-2}\lesssim i^{-1}$ and $m_1'+m_2'-n'\leq -1$:
  \begin{align*}
&\sum_{b=1}^{j-i+1}\sum_{k: M_k\in A_{j+1}^{i+b}}\lambda_{k}\frac{\det(M_k)}{\det(A)}|A^{-1}-M_{k}^{-1}|\lesssim\sum_{b=1}^{j-i+1}\sum_{k: M_k\in A_{j+1}^{i+b}}\lambda_{k} j\frac{(i+b)^{m_1'}}{i^{m_1'}}\\
&\leq \sum_{b=1}^{j-i+1}\nu(A_{j+1}^{i+b})j\frac{(i+b)^{m_1'}}{i^{m_1'}}= \sum_{b=1}^{j-i+1}j^{2(m_1'+m_2'-n')+1}  \frac{i^{m_1'}}{(i+b)^{m_1'+2} }\frac{(i+b)^{m_1'}}{i^{m_1'}}\lesssim j^{-1}i^{-1}\leq 1.
\end{align*}
Using that $i\leq j$ and $m_1'+m_2'\leq n'-1$, we have
\begin{align*}
&\sum_{k: M_k\in B_{j+1}^{j+1}}\lambda_{k}\frac{\det(M_k)}{\det(A)}|A^{-1}-M_{k}^{-1}|\lesssim \sum_{k: M_k\in B_{j+1}^{j+1}}\lambda_{k}j^{2n'-2m_2'-m_1'+1}i^{-m_1'}\\
&=\nu(B_{j+1}^{j+1})j^{2n'-2m_2'-m_1'+1}i^{-m_1'}\lesssim j^{2m_1'+3m_2'-3n'}i^{m_1'}j^{2n'-2m_2'-m_1'+1}i^{-m_1'}= j^{m_1'+m_2'-n'+1}\leq 1.
\end{align*}
 In the same way as before we bound the sum over the sets $S_{j+1,i}^{a}$ and $S_{j+1,i}^{a+b}$; using that in this case we have $m_1'=m_1$ and $n'=n$,  we get
  \begin{align*}
&\sum_{a=m_2+1}^{n-m_1-1}\sum_{k: M_k\in S_{j+1,i}^{a}}\lambda_{k}\frac{\det(M_k)}{\det(A)}|A^{-1}-M_{k}^{-1}|\lesssim\sum_{a=m_2+1}^{n-m_1-1}\sum_{k: M_k\in S_{j+1,i}^{a}}\lambda_{k} j^{n'-m_1'-a+1}i^{n'-a-m_1'}\\
&= \sum_{a=m_2+1}^{n-m_1-1}\nu(S_{j+1,i}^{a})j^{n'-m_1'-a+1}i^{n'-a-m_1'}\lesssim \sum_{a=m_2+1}^{n-m_1-1}j^{2(m_1'+a-n')} i ^{m_1'+a-n'}j^{n'-m_1'-a+1}i^{n'-a-m_1'}\\
&=\sum_{a=m_2+1}^{n-m_1-1}j^{m_1'-n'+a+1}\lesssim 1
\end{align*}
 and
  \begin{align*}
&\sum_{b=1}^{j-i+1}\sum_{a=m_2+1}^{n-m_1-1}\sum_{k: M_k\in S_{j+1,i+b}^{a}}\lambda_{k}\frac{\det(M_k)}{\det(A)}|A^{-1}-M_{k}^{-1}|\\
&\lesssim\sum_{b=1}^{j-i+1}\sum_{a=m_2+1}^{n-m_1-1}\sum_{k: M_k\in S_{j+1,i+b}^{a}}\lambda_{k}j^{n'-m_1'-a+1} (i+b)^{n'-a}i^{-m_1'}\\
&= \sum_{b=1}^{j-i+1}\sum_{a=m_2+1}^{n-m_1-1}\nu(S_{j+1,i+b}^{a})j^{n'-m_1'-a+1} (i+b)^{n'-a}i^{-m_1'}\\
&\lesssim \sum_{b=1}^{j-i+1}\sum_{a=m_2+1}^{n-m_1-1}j^{m_1'+2m_2'-n'-a+1} (i+b)^{2m_2'-2a} \lesssim j^{-1}  \sum_{b=1}^{j-i+1}(i+b)^{-2}\lesssim 1.
\end{align*}
Therefore, we get 
\[\sum_{k=1}^{N}\lambda_{k}\frac{\det(M_k)}{\det(A)}|A^{-1}-M_{k}^{-1}|\lesssim 1.\]
The proof is finished.
\end{proof}
The proof of the next lemma is analogous to the one of Lemma \ref{lemma laminado Aji}, but using Lemma \ref{lemma laminado Aki} in the induction instead of Lemma \ref{lemma laminado Bki}.
\begin{lem}\label{lemma laminado Bkj}
Let $i,j\in\N$, $i\leq j$, and $A\in B_i^{j}$. Then, there exists a laminate $\nu=\sum_{k=1}^{N}\lambda_{k}\delta_{M_{k}}$ such that 
\begin{enumerate}[a)]
\item $\overline{\nu}=A$,
\item $\supp (\nu)\subset A_{j+1}^{j+1}\cup\bigcup_{b=0}^{j-i+1} B_{i+b}^{j+1}\cup \bigcup_{b=0}^{j-i+1}\bigcup_{a=m_2+1}^{n-m_1-1}S_{i+b,j+1}^{a}$,
\item $\nu(A_{j+1}^{j+1})\lesssim  j^{m_1'+m_2'-n'}$, 
\item $\nu(B_{i}^{j+1})- \left(\frac{j+2}{j+3}\right)^{n'-m_2'}\lesssim(ij)^{-2}$,
\item $\nu(B_{i+b}^{j+1})\lesssim \left((i+b)(j+3)\right)^{2(m_1'+m_2'-n')} $, for $b\in\{1,\ldots,j-i+1\}$,
\item  $\nu(S_{i,j+1}^{a})\lesssim j^{m_2'-a}$ for $a\in\{m_2+1,\ldots,n-m_1-1\}$,
\item  $\nu(S_{i+b,j+1}^{a})\lesssim j^{m_1'+m_2'-n'}\left((i+b-1)^{2}(j+1)\right)^{m_1'+a-n'}$ for $a\in\{m_2+1,\ldots,n-m_1-1\}$ and $b\in\{1,\ldots,j-i+1\}$,
\item  $\sum_{k=1}^{N}\lambda_{k}|A-M_{k}|\lesssim 1$,
\item $\frac{1}{\det (A)}\sum_{k=1}^{N}\lambda_{k}\det (M_{k})|A^{-1}-M_{k}^{-1}|\lesssim j^{-2}$.
\end{enumerate}
\end{lem}

We recall now two results from \cite{FaMoOl16} that will be used to prove Theorem \ref{theorem existence}.

\begin{prop}\label{prop: laminate-bisobolev}

Let $N\in \N$, $A_1,\ldots, A_N\in \Gamma_+$ and $L\geq 1$ be such that
\[\sigma_{1}(A_i)\geq L^{-1} , \qquad |A_i|\leq L , \qquad i=1,\ldots, N.\]
Consider $\alpha_1, \ldots, \alpha_N \geq 0$ such that $\nu := \sum_{i=1}^{N}\alpha_{i}\delta_{A_{i}}$ is in $\mathcal{L}(\R^{n\times n})$ and call $A:=\overline{\nu}$.
Then, for every $\alpha\in (0,1)$, $0<\delta<\frac{1}{2}\min\{L^{-1},\min_{1\leq i<j\leq N}|A_{i}-A_{j}|\}$ and every bounded open set $\Omega\subset\R^{n}$, there exists a piecewise affine bi-Lipschitz homeomorphism $f:\Omega\to A\Om$ such that

\begin{enumerate}[(a)]
\item\label{ f is gradient} $f=\nabla u$ for some $u\in W^{2,\infty}(\Omega)$,
\item\label{frontera} $f(x)=Ax$ for $x\in\partial\Omega$,
\item\label{alpha f} $\|f-A\|_{C^{\alpha}(\overline{\Omega})}<\delta$,
\item\label{alpha f inv} $\|f^{-1}-A^{-1}\|_{C^{\alpha}(A\overline{\Omega})}<\delta$,
\item\label{gradient f} $|\{x\in\Omega: | D f(x)-A_{i}|<\delta \}|=\alpha_{i}|\Omega|$ for all $i=1,\ldots,N$.
\end{enumerate} 
\end{prop}
 \begin{lem}\label{lemma glue homeomorphisms}
Let $f:\overline{\Om}\to \R^n$ be a homeomorphism such that $f$ and $f^{-1}$ are $C^{\alpha}$ for some $\alpha \in (0,1]$.
Let $\{ \omega_{i} \}_{i \in \N}\subset\Om$ be a family of pairwise disjoint open sets, and for each $i\in \N$ let $g_i:\overline{\omega_i}\to f(\overline{\omega_i})$ be a homeomorphism such that $g_i=f$ on $\partial \omega_i$,
 \[\sup_{i\in\N}\|f-g_i\|_{C^{\alpha}(\overline{\omega_i})} <\infty  \quad \text{and} \quad \sup_{i\in\N}\|f^{-1}-g_i^{-1}\|_{C^{\alpha}(f(\overline{\omega_i}))} <\infty .\]
Then, the function
\[\tilde{f} (x) :=\begin{cases}
f(x) &\text{if } x \in \overline{\Om}\setminus\bigcup_{i\in\N}\omega_i,\\
g_i (x) &\text{if } x \in \omega_i \text{ for some } i \in \N
\end{cases}\]
 is a homeomorphism between $\overline{\Om}$ and $f(\overline{\Om})$ such that $\tilde{f}$ and $\tilde{f}^{-1}$ are $C^{\alpha}$ and
\[
  \| f- \tilde{f}\|_{C^{\alpha} (\overline{\Omega})} \leq 2 \sup_{i\in\N}\|f-g_i\|_{C^{\alpha}(\overline{\omega_i})} , \qquad
 \| f^{-1} - \tilde{f}^{-1}\|_{C^{\alpha} (\overline{\Omega})} \leq 2 \sup_{i\in\N}\|f^{-1}-g_i^{-1}\|_{C^{\alpha}(f(\overline{\omega_i}))} .
\]
\end{lem}

We construct now two families of uniformly bounded constants that will be used in the proof of Theorem \ref{theorem existence}.

\begin{lem}\label{lemma constants}
Let $\ve'>0$ and $\tilde{C}>1$. Define the sequences $\{C_{j,i}^{1}\}_{\substack{j\in\N\\i=1,\ldots,j}}$, $\{C_{j,i}^{2}\}_{\substack{j\in\N\\i=0,\ldots,j}}$ and $\{M_{j}\}_{j\in\N}$ as follows:
\begin{enumerate}[a)]
\item $C_{1,1}^{1}=4^{n}$, $C_{1,1}^{2}=4^{n}$, $C_{j,0}^{2}=0$ for $j\in\N$.

Given $j\in\N$, assuming $C_{j,i}^{1}$ and $C_{j,i}^{2}$ have been defined for all $i\in\{1,\ldots,j\}$, set
\item  $M_{j}=\max_{i=1,\ldots,j}C_{j,i}^{1}$,
\item\label{lem Cj,i1} $C_{j+1,i}^{1}=C_{j,i}^{1}+j^{-2}\tilde{C}M_{j}+\tilde{C}C_{j,i-1}^{2}(j+2-i)^{-2},$ for $i=1,\ldots,j$,
\item\label{lem Cj,j1} $C_{j+1,j+1}^{1}=\tilde{C}\left(M_{j}+C_{j,j}^{2}\right)j^{-\ve'},$
\item\label{lem Cj,i2} $C_{j+1,i}^{2}=C_{j,i-1}^{2}\left(1+\tilde{C}i^{-2}\right)+\tilde{C}M_{j}i^{-2+\ve'},$ for $i=1,\ldots,j$,
\item\label{lem Cj,j2} $C_{j+1,j+1}^{2}=C_{j,j}^{2}\left(1+\tilde{C}j^{-2}\right)+\tilde{C}M_{j}j^{-4},$
\end{enumerate}
Then
\[\sup_{j\in\N}M_{j}<\infty \text{ and }\sup_{\substack{j\in\N\\i=0,\ldots,j}}C_{j,i}^{2}<\infty.\]
\end{lem}
\begin{proof}
Clearly $M_{j}\leq M_{j+1}$ for $j\in\N$. From \emph{\ref{lem Cj,j2})} we obtain by induction on $j$ that
\[C_{j+1,j+1}^{2}=C_{1,1}^{2}\prod_{k=1}^{j}(1+\tilde{C}k^{-2})+\sum_{\ell=1}^{j}\left[\prod_{k=\ell+1}^{j}(1+\tilde{C}k^{-2})\right]\tilde{C}M_{\ell}\ell^{-4}.\]
Hence
\begin{equation}\label{const jj2}
C_{j+1,j+1}^{2}\lesssim \sum_{\ell=1}^{j}M_{\ell}\ell^{-4}\lesssim M_{j}.
\end{equation}
Therefore, from \emph{\ref{lem Cj,j1})} we get
\begin{equation}\label{const jj1}
C_{j+1,j+1}^{1}\lesssim M_{j}j^{-\ve'}.
\end{equation}
On the other hand, from \emph{\ref{lem Cj,i2})} by induction on $j$, we have, for $i=1,\ldots,j$, that
\begin{equation}\label{const ji2}
C_{j+1,i}^{2}=\sum_{\ell=1}^{i}\left[\prod_{k=\ell+1}^{i}(1+\tilde{C}k^{-2})\right]\tilde{C}M_{j+\ell-i}\ell^{-2+\ve'}\lesssim \sum_{\ell=1}^{i}M_{j+\ell-i}\ell^{-2+\ve'}.
\end{equation}
By \emph{\ref{lem Cj,i1})} we obtain by induction on $j\geq i-1$ that
\begin{equation}\label{const ji1}
C_{j+1,i}^{1}=C_{i,i}^{1}+\sum_{\ell=i}^{j}\tilde{C}\left(M_{\ell}\ell^{-2}+(\ell+2-i)^{-2}C_{\ell,i-1}^{2}\right).
\end{equation}
For $r\geq 1$ we use that
\[-s^{2}+(r-2)s+2r\geq\begin{cases}
\frac{r(s+1)}{2}&\text{for }0\leq s\leq \frac{r}{2},\\
\frac{r(r-s)}{2}&\text{for }\frac{r}{2}\leq s\leq r-1
\end{cases}\]
to get
\begin{equation}\label{const sum}
\begin{split}
&\sum_{\ell=i}^{j}\left[(\ell+2-i)^{-2}\sum_{k=1}^{i-1}M_{\ell+k-i}k^{-2+\ve'}\right]\leq \sum_{r=1}^{j-1}M_{r}\sum_{s=0}^{r-1}(s+2)^{-2}(r-s)^{-2+\ve'}\\
&\leq \sum_{r=1}^{j-1}M_{r}\sum_{s=0}^{r-1}(-s^{2}+(r-2)s+2r)^{-2+\ve'}\lesssim \sum_{\ell=1}^{j-1}M_{\ell}\ell^{-2+\ve'}.
\end{split}
\end{equation}
Now, we use (\ref{const jj1}), (\ref{const ji2}), (\ref{const ji1}) and (\ref{const sum}) to get that there exists a constant $C'$ depending only on $n$ such that
\begin{equation*}
C_{j+1,i}^{1}\leq C_{i,i}^{1}+ C'\sum_{\ell=1}^{j}M_{\ell}\ell^{-2+\ve'}.
\end{equation*}
Let $j_{\ve'}$ be such that for all $j\geq j_{\ve'}$
\[C_{j+1,j+1}^{1}\leq M_{j},\]
which is possible thanks to (\ref{const jj1}). Define the family of constants $\{C_{j,i}^{3}\}_{\substack{j\in\N\\i=1,\ldots,j}}$ as follows
\begin{equation*}
C_{j,i}^{3}=
\begin{cases}
C_{j,i}^{1}&\text{if }j\leq j_{\ve'}\text{ or }i=j,\\
C_{i,i}^{1}+C'\sum_{\ell=1}^{j-1}M_{\ell}\ell^{-2+\ve'}&\text{if }j>\max\{i,j_{\ve'}\}.
\end{cases}
\end{equation*}
Hence, $C_{j,i}^{1}\leq C_{j,i}^{3}\text{ for all }j\in\N\text{ and }i=1,\ldots,j$. Define $M_{j}'=\max_{i=1,\ldots,j}C_{j,i}^{3}$ and let $i_{j}\in\{1,\ldots,j\}$ be such that $M_{j}'=C_{j,i_j}^{3}$. First, we note that $M_{j}'<M_{j+1}'$ for all $j>j_{\ve'}$. Fix $j>j_{\ve'}$; as
\[C_{j+1,j+1}^{3}=C_{j+1,j+1}^{1}\leq M_{j}\leq M_{j}'\text{ and } M_{j}'<M_{j+1}',\]
it is clear that $i_{j+1}\leq j$. We also see that
\[C_{j+1,i_{j}}^{3}=C_{j,i_{j}}^{3}+C'M_{j}j^{-2+\ve'}=M_{j}'+C'M_{j}j^{-2+\ve'},\]
and for all $i\leq j$ we have
\[C_{j+1,i}^{3}=C_{j,i}^{3}+C'M_{j}j^{-2+\ve'}. \]
Therefore $M_{j+1}'=C_{j+1,i_{j}}^{3}$, so we can take $i_{j+1}=i_{j}$. By induction, we can take $i_{j}=i_{j_{\ve'}+1}$ for all $j>j_{\ve'}$, and $M_{j}'=C_{j,i_{j_{\ve'}+1}}^{3}$. Hence, there exists $C>0$ such that
\[M_{j+1}'\leq C \sum_{\ell=1}^{j}M_{\ell}\ell^{-2+\ve'}\leq C\sum_{\ell=1}^{j}M_{\ell}'\ell^{-2+\ve'}.\]
Define for $j\leq j_{\ve'}+1$,  $\tilde{M}_{j}=M_{j}'$ and for $j\geq j_{\ve'}+1$
\[\tilde{M}_{j+1}=C \sum_{\ell=1}^{j}\tilde{M}_{\ell}\ell^{-2+\ve'},\]
so $M_{j}'\leq \tilde{M}_{j}$ for all $j\in\N$. Now, we observe that for $j\geq j_{\ve'}+2$ we have 
\[\tilde{M}_{j+1}=\tilde{M}_{j}(1+Cj^{-2+\ve'}).\]
Therefore 
\[\sup_{j\in\N}M_{j}\leq \sup_{j\in\N}M_{j}'\leq \sup_{j\in\N}\tilde{M}_{j}<\infty.\]
Using (\ref{const jj2}), (\ref{const ji2}) and that $\sup_{j\in \N}M_{j}<\infty$  we have 
\[\sup_{\substack{j\in\N\\i=1,\ldots,j}}C_{j,i}^{2}<\infty.\]
Here concludes the proof.
\end{proof}
Finally, we combine all the previous results to prove Theorem \ref{theorem existence}.
\begin{proof}[Proof of Theorem \ref{theorem existence}]
In this proof, expressions like $\{x\in\Omega: f(x)\in A\}$ will be abbreviated as   $\{f(x)\in A\}$. Given $\ve'>0$ small enough to have $\sum_{k=1}^{\infty}k^{-2+\ve'}<2$, we will construct a sequence $\{f_j\}_{j\in\N}\subset W^{1,1}\left(\Omega,\Omega\right)$ of piecewise affine Lipschitz homeomorphisms such that  $f_0=\id$ and, when we take $\{C_{j,i}^{1}\}_{\substack{j\in\N\\i=1,\ldots,j}}$, $\{C_{j,i}^{2}\}_{\substack{j\in\N\\i=0,\ldots,j}}$ the families of constants in Lemma \ref{lemma constants} and we denote $\Omega_{S}^{j}=\{Df_j(x)\in S\}$ for each $S\subset\Gamma_{+}$, we have 
\begin{enumerate}[i)]
\item\label{pro th ex border} $f_j=\id$ on $\partial\Omega$,
\item\label{pro th ex holder} $\|f_{j}-f_{j-1}\|_{C^{\alpha}(\overline{\Omega})}<2^{-j}\ve$ and  $\|f_{j}^{-1}-f_{j-1}^{-1}\|_{C^{\alpha}(\overline{\Omega})}<2^{-j}\ve$,
\item\label{pro th ex sup} $Df_j(x)\in \bigcup_{i=1}^{j}\bigcup_{a=m_2+1}^{n-m_1-1}\left( A_{j}^{i}\cup B_{i}^{j}\cup S_{j,i}^{a}\cup S_{i,j}^{a} \right)$,
\item\label{pro th ex sob conv} $\int_{\Omega}|Df_{j}(x)-Df_{j-1}(x)|dx\lesssim j^{-2}|\Omega|$,
\item\label{pro th ex sob inv conv} $\int_{\Omega}|Df_{j}^{-1}(y)-Df_{j-1}^{-1}(y)|dy\lesssim j^{-2}|\Omega|$,
\item\label{pro th ex A} $\frac{|\Omega_{A_{j}^{i}}^{j}|}{|\Omega|}\leq C_{j,i}^{1} i^{-m_1'-2+\ve'}$ for $i=1,\ldots,j$,
\item\label{pro th ex B} $\frac{|\Omega_{B_{i}^{j}}^{j}|}{|\Omega|}\leq C_{j,i}^{1} i^{-2+\ve'}(j+2)^{m_2'-n'}$ for $i=1,\ldots,j$,
\item\label{pro th ex Sji} $\frac{|\Omega_{S_{j,i}^{a}}^{j}|}{|\Omega|}\leq C_{j,i}^{2} (i+2)^{a-n'}(j+1-i)^{-2}$ for $i=1,\ldots,j$ and $a\in\{m_2+1,\ldots, n-m_1-1\}$,
\item\label{pro th ex Sij} $\frac{|\Omega_{ S_{i,j}^{a}}^{j}|}{|\Omega|}\leq C_{j,i}^{2} (j+2)^{a-n'}(j+1-i)^{-2}$ for $i=1,\ldots,j$ and $a\in\{m_2+1,\ldots, n-m_1-1\}$.
\end{enumerate}
One constructed such sequence $\{f_j\}$, we have that it converges in the $C^{\alpha}$ and in the $W^{1,1}$ norm to a bi-Sobolev homeomorphism $f:\Omega\to\Omega$; see, if necessary the proof of \cite[Th.\,1]{FaMoOl16} for the details of the limit passage. Moreover, it is immediate from \ref{pro th ex border}), \ref{pro th ex holder}), that $f$ satisfies \emph{\ref{th ex border})} and \emph{\ref{th ex holder})} of Theorem \ref{theorem existence}. Recall that the bounds of $S_{i,k}^{a}$ only have sense if $m_1+m_2\leq n-1$; otherwise, these sets are empty. From \ref{pro th ex sup}), \ref{pro th ex B}), \ref{pro th ex Sji}) and \ref{pro th ex Sij}) we obtain
\begin{equation}\label{pro th lim rank}
1-\frac{|\Omega_{ \bigcup_{i=1}^{j} A_{j}^{i}}^{j}|}{|\Omega|}\lesssim \sum_{i=1}^{j}\left[i^{-2+\ve'}j^{m_2'-n'}+(j+1-i)^{-2}\left(\sum_{a=m_2+1}^{n-m_1-1}i^{a-n'}+j^{a-n'}\right)\right]\lesssim j^{-1}.
\end{equation}
For a subsequence, $Df_j \to D f$ a.e., so, thanks to the continuity of the singular values and using that $\Gamma_{+}$ is closed we obtain that $Df \in \Gamma_+$ and, thanks to \eqref{pro th lim rank}, we also get that $\rank \left( D f \right)=m_1$ a.e.\ in $\Om$. On the other hand,
\begin{align*}
&1-\frac{\left|\left\{ Df_{j}^{-1}(y)\in \bigcup_{i=1}^{j}\left( B_{i}^{j}\right)^{-1}\right\}\right|}{|\Omega|}\\
&=\frac{\left|\left\{Df_{j}^{-1}(y)\in \bigcup_{i=1}^{j}\bigcup_{a=m_2+1}^{n-m_1-1}\left( A_{j}^{i}\cup S_{i,j}^{a}\cup S_{j,i}^{a}\right)^{-1}\right\}\right|}{|\Omega|}\\
&=\frac{1}{|\Omega|}\left|f_{j}\left(\left\{ Df_j(x)\in \bigcup_{i=1}^{j} \bigcup_{a=m_2+1}^{n-m_1-1} A_{j}^{i}\cup S_{i,j}^{a}\cup S_{j,i}^{a}\right\}\right)\right|\\
&=\frac{1}{|\Omega|}\int_{\left\{Df_j(x)\in \bigcup_{i=1}^{j} \bigcup_{a=m_2+1}^{n-m_1-1} A_{j}^{i}\cup S_{i,j}^{a}\cup S_{j,i}^{a}\right\}}\det(Df_{j}(x))dx.
\end{align*}
Now, we split the integral over the different sets and we use the control that we have over the determinant in the different sets (the second part with $M\in B_{i}^{j}$ will be used later):
\begin{equation}\label{pro th ex control det}
\det(M)\lesssim
\begin{cases}
j^{m_1'-n'}i^{m_1'}&\text{if }M\in A_{j}^{i},\\
j^{n'-m_2'}i^{-m_2'}&\text{if }M\in B_{i}^{j},\\
k^{-a}i^{n'-a}&\text{if }M\in S_{k,i}^{a}.
\end{cases}
\end{equation}
Therefore using  \ref{pro th ex A}), \ref{pro th ex Sji}) and \ref{pro th ex Sij}) we get
\begin{align*}
&1-\frac{\left|\left\{Df_{j}^{-1}(x)\in \bigcup_{i=1}^{j}\left( B_{i}^{j}\right)^{-1}\right\}\right|}{|\Omega|}\\
&\lesssim \sum_{i=1}^{j}\Bigg[i^{-m_1'-2+\ve'}j^{m_1'-n'}i^{m_1'}+(j+1-i)^{-2}\sum_{a=m_2+1}^{n-m_1-1}\left(i^{a-n'}j^{-a}i^{n'-a}+j^{a-n'}i^{-a}j^{n'-a}\right)\Bigg]\lesssim j^{-1}.
\end{align*}
The same argument as before shows that $\rank(Df^{-1}(y))=m_2$ a.e. $y\in\Omega$. Hence, \emph{\ref{th ex rank})} of Theorem \ref{theorem existence} is proved.

Let $\ve,t>0$ and pick $j>t$; then using that 
\begin{equation*}
|M|\leq
\begin{cases}
i+2&\text{if }M\in A_{j}^{i},\\
j+2&\text{if }M\in B_{i}^{j},\\
i+2&\text{if }M\in S_{k,i}^{a},
\end{cases}
\end{equation*}
we have
\begin{align*}
&\frac{|\{ |Df_{j}(x)|>t\}|}{|\Omega|}\leq \sum_{i=\floor{t}-1}^{j}\left[\frac{|\Omega_{ A_{j}^{i}}^{j}|}{|\Omega|}+\sum_{a=m_2+1}^{n-m_1-1}\frac{|\Omega_{ S_{j,i}^{a}}^{j}|}{|\Omega|}\right]+\sum_{i=1}^{j}\left[\frac{|\Omega_{B_{i}^{j}}^{j}|}{|\Omega|}+\sum_{a=m_2+1}^{n-m_1-1}\frac{|\Omega_{ S_{i,j}^{a}}^{j}|}{|\Omega|}\right]\\
&\lesssim \sum_{i=\floor{t}-1}^{j}\left[i^{-m_1'-2+\ve'}+\sum_{a=m_2+1}^{n-m_1-1}i^{a-n'}(j+1-i)^{-2}\right]\\
&\quad+\sum_{i=1}^{j}\left[i^{-2+\ve'}j^{m_2'-n'}+ \sum_{a=m_2+1}^{n-m_1-1}j^{a-n'}(j+1-i)^{-2}\right].
\end{align*}
Hence, using $\sum_{i=\floor{t}-1}^{j}i^{-m_1'-1}(j+1-i)^{-2}\lesssim t^{-m_1'-1}\leq t^{-m_1'-1+\ve'}$ and $m_1'+m_2'\leq n'-1$ we get 
\[\frac{|\{|Df_{j}(x)|>t\}|}{|\Omega|}\lesssim t^{-m_1'-1+\ve'},\]
and, hence, since we have proved the bound for all $\ve'>0$ we have $f\in W^{1,p}(\Omega,\R^{n})$ for all $p<m_1'+1$.

Next, for the inverse we will use the bounds
\begin{equation*}
|M^{-1}|\leq
\begin{cases}
j+2&\text{if }M\in A_{j}^{i},\\
i+2&\text{if }M\in B_{i}^{j},\\
k+2&\text{if }M\in S_{k,i}^{a}.
\end{cases}
\end{equation*}
Therefore we get
\begin{align*}
&\frac{|\{ |Df_{j}^{-1}(y)|>t\}|}{|\Omega|}=\frac{|f_{j}\left(\{|Df_{j}(x)^{-1}|>t\}\right)|}{|\Omega|}\\
&\leq\sum_{i=1}^{j}\left[\frac{|f_{j}(\Omega_{ A_{j}^{i}}^{j})|}{|\Omega|}+\sum_{a=m_2+1}^{n-m_1-1} \frac{|f_{j}(\Omega_{ S_{j,i}^{a}}^{j})|}{|\Omega|}\right]+\sum_{i=\floor{t}-1}^{j}\left[\frac{|f_{j}(\Omega_{ B_{i}^{j}}^{j})|}{|\Omega|}+\sum_{a=m_2+1}^{n-m_1-1} \frac{|f_{j}(\Omega_{S_{i,j}^{a}}^{j})|}{|\Omega|}\right].
\end{align*}
Therefore, using \eqref{pro th ex control det} we get
\begin{align*}
&\frac{|\{ |Df_{j}^{-1}(y)|>t\}|}{|\Omega|}\lesssim \sum_{i=1}^{j}\left[j^{m_1'-n'}i^{m_1'}i^{-m_1'-2+\ve'}+ \sum_{a=m_2+1}^{n-m_1-1}j^{-a}i^{n'-a}i^{a-n'}(j+1-i)^{-2}\right]\\
&+\sum_{i=\floor{t}-1}^{j}\left[ j^{n'-m_2'}i^{-m_2'}i^{-2+\ve'}j^{m_2'-n'}+ \sum_{a=m_2+1}^{n-m_1-1}i^{-a}j^{n'-a}j^{a-n'}(j+1-i)^{-2}\right]\lesssim t^{-m_2'-1+\ve'}.
\end{align*}
So, we have $f^{-1}\in W^{1,q}(\Omega,\R^{n})$ for all $q<m_2+1$, and therefore, part \emph{\ref{th ex int}}) of Theorem \ref{theorem existence} is proved.

Hence, to prove the theorem it is enough to construct the sequence $\{f_j\}$. Let $f_0=\id$ and proceeding as in Lemma \ref{lemma laminado Aki} we construct a laminate $\nu_1$ such that $\overline{\nu}_1=I$ and 
\[\supp(\nu_1)\subset A_{1}^{1}\cup B_{1}^{1}\cup\bigcup_{a=m_2+1}^{n-m_1-1} S_{1,1}^{a}.\]
Now apply Proposition \ref{prop: laminate-bisobolev} with $\delta$ small enough to have $Df_{1}(x)\in A_{1}^{1}\cup B_{1}^{1}\cup\bigcup_{a=m_2+1}^{n-m_1-1} S_{1,1}^{a}$ for almost every $x\in\Omega$; this is possible because Proposition \ref{prop: laminate-bisobolev} gives us that $Df(x)\in\Gamma_{+}$ a.e. $x\in\Omega$; and, for $i,k\in\N$, $a\in\{m_2+1,\ldots,n-m_1-1\}$, the sets $A_{i}^{k},B_{i}^{k},S_{i,k}^{a}$ are open in $\Gamma_{+}$. We do now the inductive step: suppose that we have constructed $f_j$. As $f_j$ is piecewise affine, there exist  families $\{A_k\}_{k\in\N}\subset\R^{n\times n}$, $\{b_{k}\}_{k\in\N}\subset\R^{n}$ and $\{\omega_{k}\}_{k\in\N}$ of disjoint sets such that $|\Omega\setminus \left(\bigcup_{k\in\N}\omega_{k}\right)|=0$ and 
\[f_j(x)=A_{k}x+b_{k}\text{ in }\omega_{k}.\]
For each $k$, let $\nu_k$ the laminate given by Lemma \ref{lemma laminado Ski} if $A_{k}\in \bigcup_{i=1}^{j}\bigcup_{a=m_2+1}^{n-m_1-1} S_{j,i}^{a}\cup S_{i,j}^{a}$, the one given by Lemma \ref{lemma laminado Aji} if $A_{k}\in \bigcup_{i=1}^{j}A_{j}^{i}$ and the one given by Lemma \ref{lemma laminado Bkj} if $A_{k}\in \bigcup_{i=1}^{j}B_{i}^{j}$. For each $k$, let $g_k$ be the homeomorphism given by Proposition \ref{prop: laminate-bisobolev} corresponding to $\nu_k$ that is equal to $A_k x+b_k$ on the border of $\omega_k$ with $\delta_k>0$ being as small as we need in the rest of the proof. First let $\delta_k$ be small enough to get 
\[Dg_k(x)\in \bigcup_{i=1}^{j+1}\bigcup_{a=m_2+1}^{n-m_1-1}\left( A_{j+1}^{i}\cup B_{i}^{j+1}\cup S_{j+1,i}^{a}\cup S_{i,j+1}^{a}\right) \text{ for almost every }x\in\omega_k,\]
and such that for all $k\in\N$ we have $\delta_k<2^{-j-1}\ve$. Then, we define
\[f_{j+1}=
\begin{cases}
g_k&\text{in }\omega_k\text{ for some }k\in\N,\\
f_j&\text{in }\Omega\setminus\bigcup_{k\in\N}\omega_k.
\end{cases}\]
It is obvious that $f_{j+1}$ satisfy \ref{pro th ex border}) and \ref{pro th ex sup}); using Lemma \ref{lemma glue homeomorphisms} we see that it is a homeomorphism, and by Proposition \ref{prop: laminate-bisobolev} and Lemma \ref{lemma glue homeomorphisms} we have \ref{pro th ex holder}). Now we prove \ref{pro th ex sob conv}) and \ref{pro th ex sob inv conv}). Choose $\delta_k$ such that, if $\nu_k=\sum_{\ell=1}^{N_k}\lambda_{k,\ell}\delta_{M_{k,\ell}}$, 
\[\delta_k<\min_{\ell\in\{1,\ldots,N_k\}}|A_k-M_{k,\ell}|,\]
and for $\ell=1,\ldots,N_k$, using part \emph{(\ref{gradient f})} of Proposition \ref{prop: laminate-bisobolev} we have that it also holds
\[|Dg_{k}^{-1}(y)-M_{k,\ell}^{-1}|<|A_{k}^{-1}y-M_{k,\ell}^{-1}|\quad\text{in }g_{k}(\{x\in \omega_k: |Dg_k(x)-M_{k,\ell}|<\delta_k\}).\]
Denote by $\omega_{k,\ell}$ the set $\{x\in\omega_{k}: |Dg_{k}(x)-M_{k,\ell}|<\delta_{k}\}$. Recall, from Proposition \ref{prop: laminate-bisobolev} \emph{(\ref{gradient f})}, that $|\omega_{k,\ell}|=\lambda_{k,\ell}|\omega_k|$. Then, using parts \emph{\ref{Aji conv f})} and \emph{\ref{Aji conv f-1})} of Lemma \ref{lemma laminado Aji}, we have that for those $k\in\N$ such that $A_k\in \bigcup_{i=1}^{j}A_{j}^{i}$,
\begin{align*}
\int_{\omega_k}|A_k-Dg_k(x)|dx&\leq \sum_{\ell=1}^{N_k}\int_{\omega_{k,\ell}}\left(|A_k-M_{k,\ell}|+|M_{k,\ell}-Dg_k(x)|\right)dx\\
&\lesssim  \sum_{\ell=1}^{N_k}\int_{\omega_{k,\ell}}|A_k-M_{k,\ell}|dx=\sum_{\ell=1}^{N_k}\lambda_{k,\ell}|A_k-M_{k,\ell}||\omega_{k}|\lesssim j^{-2}|\omega_{k}|,
\end{align*}
and, also using $|g_k(\omega_{k,\ell})|\lesssim \det(M_{k,\ell})|\omega_{k,\ell}|=\det(M_{k,\ell})\lambda_{k,\ell}|\omega_{k}|$ we get
\begin{align*}
&\int_{g_{k}(\omega_k)}|A_{k}^{-1}-Dg_{k}^{-1}(y)|dx\leq \sum_{\ell=1}^{N_k}\int_{g_k(\omega_{k,\ell})}\left(|A_{k}^{-1}-M_{k,\ell}^{-1}|+|M_{k,\ell}^{-1}-Dg_{k}^{-1}(y)|\right) dy\\
&\lesssim  \sum_{\ell=1}^{N_k}\int_{g_k(\omega_{k,\ell})}|A_{k}^{-1}-M_{k,\ell}^{-1}|dx\lesssim  \sum_{\ell=1}^{N_k}\lambda_{k,\ell}\det(M_{k,\ell})|\omega_{k}||A_{k}^{-1}-M_{k,\ell}^{-1}|\lesssim \det(A_{k})|\omega_{k}|=|g_{k}(\omega_{k})|.
\end{align*}
Proceeding in the same way as before we obtain that for those $k\in \N$ such that $A_{k}\in \bigcup_{i=1}^{j}B_{i}^{j}$ we have
\[\int_{\omega_k}|A_k-Dg_k(x)|dx\lesssim |\omega_{k}|\]
and
\[\int_{g_{k}(\omega_k)}|A_{k}^{-1}-Dg_{k}^{-1}(y)|dy\lesssim j^{-2}|g_{k}(\omega_{k})|.\]
Similarly, for $k\in\N$ such that $A_{k}\in \bigcup_{i=1}^{j}\bigcup_{a=m_2+1}^{n-m_1-1}S_{j,i}^{a}\cup S_{i,j}^{a}$ we have
\[\int_{\omega_k}|A_k-Dg_k(x)|dx\lesssim |\omega_{k}|\]
and
\[\int_{g_{k}(\omega_k)}|A_{k}^{-1}-Dg_{k}^{-1}(y)|dy\lesssim |g_{k}(\omega_{k})|.\]
Now we combine the last equations, and we use that $Df_{j+1}(x)=Dg_{k}(x)$  and $Df_{j}(x)=A_{k}$ for $x$ in $\omega_{k}$, to prove \ref{pro th ex sob conv}) and \ref{pro th ex sob inv conv}):
\begin{align*}
&\int_{\Omega}|Df_{j+1}-Df_{j}|dx=\sum_{k\in\N}\int_{\omega_{k}}|Dg_{k}-A_{k}|dx\\
&\lesssim \sum_{\substack{k\in\N\\ A_{k}\in\bigcup_{i=1}^{j}A_{j}^{i}}}j^{-2}|\omega_{k}|+\sum_{\substack{k\in\N\\ A_{k}\in\bigcup_{i=1}^{j}B_{i}^{j}}}|\omega_{k}|+\sum_{\substack{k\in\N\\ A_{k}\in\bigcup_{i=1}^{j}\bigcup_{a=m_2+1}^{n-m_1-1}S_{j,i}^{a}\cup S_{i,j}^{a}}}|\omega_{k}|\\
&=j^{-2}|\Omega_{\bigcup_{i=1}^{j}A_{j}^{i}}^{j}|+|\Omega_{\bigcup_{i=1}^{j}B_{i}^{j}}^{j}|+|\Omega_{\bigcup_{i=1}^{j}\bigcup_{a=m_2+1}^{n-m_1-1}S_{j,i}^{a}\cup S_{i,j}^{a}}^{j}|.
\end{align*}
Now, we use \ref{pro th ex A})--\ref{pro th ex Sij}), that $m_2'-n'\leq -2$ and we recall that in the bounds for the sets $S_{k,i}^{a}$ we can suppose $n'=n$, $m_1'=m_1$ and $m_2'=m_2$, because otherwise they are empty. Therefore, we obtain
\begin{align*}
&\frac{1}{|\Omega|}\int_{\Omega}|Df_{j+1}-Df_{j}|dx\\
&\lesssim \sum_{i=1}^{j}\left( j^{-2}i^{-m_1'-2+\ve'}+i^{-2+\ve'}j^{m_2'-n'}+\sum_{a=m_2+1}^{n-m_1-1}\left[i^{a-n'}(j+1-i)^{-2}+j^{a-n'}(j+1-i)^{-2}\right]\right)\\
&\lesssim j^{-2}+\sum_{i=1}^{j}i^{-2}(j+1-i)^{-2}\lesssim j^{-2}.
\end{align*}
So \ref{pro th ex sob conv}) is proved. On the other hand, using $f_{j+1}(\omega_{k})=g_{k}(\omega_{k})=f_{j}(\omega_{k})$ we have
\begin{align*}
&\int_{\Omega}|Df_{j+1}^{-1}-Df_{j}^{-1}|dy=\sum_{k\in\N}\int_{g_{k}(\omega_{k})}|Dg_{k}^{-1}-A_{k}^{-1}|dy\\
&\lesssim \sum_{\substack{k\in\N\\ A_{k}\in\bigcup_{i=1}^{j}A_{j}^{i}}}|g_{k}(\omega_{k})|+\sum_{\substack{k\in\N\\ A_{k}\in\bigcup_{i=1}^{j}B_{i}^{j}}}j^{-2}|g_{k}(\omega_{k})|+\sum_{\substack{k\in\N\\ A_{k}\in\bigcup_{i=1}^{j}\bigcup_{a=m_2+1}^{n-m_1-1}S_{j,i}^{a}\cup S_{i,j}^{a}}}|g_{k}(\omega_{k})|\\
&=\left|f_{j}\left(\Omega_{\bigcup_{i=1}^{j}A_{j}^{i}}\right)\right|+j^{-2}\left|f_{j}\left(\Omega_{\bigcup_{i=1}^{j}B_{i}^{j}}\right)\right|+\left|f_{j}\left(\Omega_{\bigcup_{i=1}^{j}\bigcup_{a=m_2+1}^{n-m_1-1}S_{j,i}^{a}\cup S_{i,j}^{a}}\right)\right|.
\end{align*}
Using again \ref{pro th ex A})--\ref{pro th ex Sij}) and that we have a control over the determinant of $Df_{j}(x)$ when $x$ is in the different sets, see (\ref{pro th ex control det}), we obtain  
\begin{align*}
&\frac{1}{|\Omega|}\int_{\Omega}|Df_{j+1}^{-1}-Df_{j}^{-1}|dy\lesssim \sum_{i=1}^{j}\left(j^{m_1'-n'}i^{m_1'}i^{-m_1'-2+\ve'}+j^{-2}j^{n'-m_2'}i^{-m_2'}i^{-2+\ve'}j^{m_2'-n'}\right.\\
&\quad+\sum_{a=m_2+1}^{n-m_1-1}\left[\left.j^{-a}i^{n'-a}i^{a-n'}(j+1-i)^{-2}+i^{-a}j^{n'-a}j^{a-n'}(j+1-i)^{-2}\right]\right)\\
&\lesssim \left(j^{-2}+\sum_{i=1}^{j}i^{-2}(j+1-i)^{-2}\right)\lesssim j^{-2}.
\end{align*}
Hence, we have proved \ref{pro th ex sob inv conv}). Finally, we suppose that \ref{pro th ex A})--\ref{pro th ex Sij}) holds for $j$ and we prove them for $j+1$. Let $C$ be a constant bigger than the ones appearing in Lemmas \ref{lemma laminado Ski}, \ref{lemma laminado Aji} and \ref{lemma laminado Bkj}.

Let $i=1,\ldots,j$; from the construction of $\nu_k$ we can see that if $\supp(\nu_{k})\cap A_{j+1}^{i}\neq \emptyset$ then
 \[\overline{\nu}_{k}=A_{k}\in \bigcup_{l=1}^{i}A_{j}^{l}\cup\bigcup_{a=m_2+1}^{n-m_1-1} S_{j,i-1}^{a}.\]
 On the other hand, using part \emph{(\ref{gradient f})} of Proposition \ref{prop: laminate-bisobolev} we get
 \begin{align*}
 |\Omega_{A_{j+1}^{i}}^{j+1}|&=\sum_{k: A_{k}\in \bigcup_{l=1}^{i}A_{j}^{l}\cup \bigcup_{a=m_2+1}^{n-m_1-1}S_{j,i-1}^{a}}\sum_{\ell: M_{k,\ell}\in A_{j+1}^{i}}\lambda_{k,\ell}|\omega_{k}|\\
 &=\sum_{k: A_{k}\in \bigcup_{l=1}^{i}A_{j}^{l}\cup\bigcup_{a=m_2+1}^{n-m_1-1} S_{j,i-1}^{a}}\nu_{k}(A_{j+1}^{i})|\omega_{k}|\\
 &=\sum_{l=1}^{i}\sum_{k:A_k\in A_{j}^{l}}\nu_{k}(A_{j+1}^{i})|\omega_{k}|+\sum_{a=m_2+1}^{n-m_1-1}\sum_{k: A_{k}\in S_{j,i-1}^{a}}\nu_{k}(A_{j+1}^{i})|\omega_{k}|.
 \end{align*}
Now, we use the control that we have over $\nu_{k}(A_{j+1}^{i})$ given by Lemmas \ref{lemma laminado Ski} and \ref{lemma laminado Aji}, and also that given $S\subset \Gamma_{+}$ we have $\sum_{k: A_{k}\in S}|\omega_{k}|=|\Omega_{S}^{j}|$. Therefore, we obtain
 \begin{align*}
 |\Omega_{A_{j+1}^{i}}^{j+1}|&\leq\sum_{l=1}^{i-1}\sum_{k:A_k\in A_{j}^{l}}\!Cj^{2(m_1'+m_2'-n')}\frac{l^{m_1'}}{i^{m_1'+2}}|\omega_{k}|+\!\sum_{k:A_k\in A_{j}^{i}}\!|\omega_{k}|+\!\sum_{a=m_2+1}^{n-m_1-1}\sum_{k: A_{k}\in S_{j,i-1}^{a}}\!Ci^{a+m_1'-n'}|\omega_{k}|\\
 &=\sum_{l=1}^{i-1}Cj^{2(m_1'+m_2'-n')}\frac{l^{m_1'}}{i^{m_1'+2}}|\Omega_{A_{j}^{l}}^{j}|+|\Omega_{A_{j}^{i}}^{j}|+\sum_{a=m_2+1}^{n-m_1-1}Ci^{a+m_1'-n'}|\Omega_{S_{j,i-1}^{a}}^{j}|.
 \end{align*}
Let $M_j$ be as in Lemma \ref{lemma constants}. By induction we get
\begin{align*}
\frac{|\Omega_{A_{j+1}^{i}}^{j+1}|}{|\Omega|}&\leq \sum_{l=1}^{i-1}C_{j,l}^{1}Cj^{2(m_1'+m_2'-n')}\frac{l^{m_1'}}{i^{m_1'+2}}l^{-m_1'-2+\ve'}+C_{j,i}^{1}i^{-m_1'-2+\ve'}\\
&\quad+\sum_{a=m_2+1}^{n-m_1-1}C_{j,i-1}^{2}Ci^{a+m_1'-n'}(i+2)^{a-n'}(j+2-i)^{-2}\\
&\leq \left(C_{j,i}^{1}+j^{-2}2CM_{j}+nCC_{j,i-1}^{2}(j+2-i)^{-2}\right)i^{-m_1'-2+\ve'}\leq C_{j+1,i}^{1}i^{-m_1'-2+\ve'}.
\end{align*}
The last estimate comes from \emph{\ref{lem Cj,i1})} of Lemma \ref{lemma constants}; we will use this lemma several times in the rest of the proof with $\tilde{C}=16^{n}C$. Now, for each $i=1,\ldots, j$, and $k\in\N$ such that  $\supp(\nu_{k})\cap B_{i}^{j+1}\neq \emptyset$, we have
 \[\overline{\nu}_{k}=A_{k}\in \bigcup_{l=1}^{i}B_{l}^{j}\cup\bigcup_{a=m_2+1}^{n-m_1-1} S_{i-1,j}^{a},\]
 hence, proceeding as before we get
  \begin{align*}
\frac{|\Omega_{B_{i}^{j+1}}^{j+1}|}{|\Omega|}&\leq \sum_{l=1}^{i-1}C_{j,l}^{1}C(i(j+3))^{2(m_1'+m_2'-n')}l^{-2+\ve'}(j+2)^{m_2'-n'}\\
&\quad+C_{j,i}^{1}\left(\left(\frac{j+2}{j+3}\right)^{n'-m_2'}+C(ij)^{-2}\right)i^{-2+\ve'}(j+2)^{m_2'-n'}\\
&\quad+\sum_{a=m_2+1}^{n-m_1-1}C_{j,i-1}^{2}C(i^{2}j)^{m_2'-a}(j+2)^{a-n'}(j+2-i)^{-2}\\
&\leq \left(C_{j,i}^{1}+j^{-2}2^{n+1}CM_{j}+n2^{2n}CC_{j,i-1}^{2}(j+2-i)^{-2}\right)i^{-2+\ve'}(j+3)^{m_2'-n'}\\
&\leq C_{j+1,i}^{1}i^{-2+\ve'}(j+3)^{m_2'-n'}.
\end{align*}
Now, we use that those $k\in\N$ such that $\nu_k$ satisfy
\[\supp(\nu_{k})\cap A_{j+1}^{j+1}\neq \emptyset\text{ or }\supp(\nu_{k})\cap B_{j+1}^{j+1}\neq \emptyset\]
are those that satisfy
\[\overline{\nu}_{k}=A_{k}\in \bigcup_{l=1}^{j}\left(A_{j}^{l}\cup B_{l}^{j}\right)\cup\bigcup_{a=m_2+1}^{n-m_1-1} S_{j,j}^{a}.\]
Therefore, we obtain
\begin{align*}
\frac{|\Omega_{A_{j+1}^{j+1}}^{j+1}|}{|\Omega|}&\leq \sum_{l=1}^{j} C_{j,l}^{1}C\left(j^{2(m_1'+m_2'-n')}\frac{l^{m_1'}}{j^{m_1'+2}}l^{-m_1'-2+\ve'}+j^{m_1'+m_2'-n'}l^{-2+\ve'}(j+2)^{m_2'-n'}\right)\\
&\quad+\sum_{a=m_2+1}^{n-m_1-1}C_{j,j}^{2}Cj^{a+m_1'-n'}(j+2)^{a-n'}\leq M_{j}4Cj^{-m_1'-2}+nCC_{j,j}^{2}j^{-m_1'-2}\\
&= \left(4CM_{j}+nCC_{j,j}^{2}\right)j^{-m_1'-2}\leq C_{j+1,j+1}^{1}(j+1)^{-m_1'-2+\ve'}
\end{align*}
and
\begin{align*}
\frac{|\Omega_{B_{j+1}^{j+1}}^{j+1}|}{|\Omega|}&\leq \sum_{l=1}^{j} C_{j,l}^{1}C\left(j^{2m_1'+3m_2'-3n'}l^{m_1'}l^{-m_1'-2+\ve'}+(j(j+1))^{2(m_1'+m_2'-n')}l^{-2+\ve}(j+2)^{m_2'-n'}\right)\\
&\quad+\sum_{a=m_2+1}^{n-m_1-1}C_{j,j}^{2}Cj^{3(m_2'-a)}(j+2)^{a-n'}\leq M_{j}4Cj^{m_2'-n'-2}+nCC_{j,j}^{2}j^{m_2'-n'-2}\\
&= \left(4CM_{j}+nCC_{j,j}^{2}\right)j^{m_2'-n'-2}\leq C_{j+1,j+1}^{1}(j+1)^{-2+\ve'}(j+3)^{m_2'-n'}.
\end{align*}
Hence, \ref{pro th ex A}) and \ref{pro th ex B}) are proved for $j+1$.

Finally, proceeding as before, we only have to prove \ref{pro th ex Sji}) and \ref{pro th ex Sij}) for $j+1$. For each $i=1,\ldots, j$, $a=m_2+1,\ldots,n-m_1-1$  and $k\in\N$ such that  $\supp(\nu_{k})\cap S_{j+1,i}^{a}\neq \emptyset$, we have
 \[\overline{\nu}_{k}=A_{k}\in \bigcup_{l=1}^{i}A_{j}^{l}\cup\bigcup_{a=m_2+1}^{n-m_1-1} S_{j,i-1}^{a},\]
 hence, using $2m_2-a-n\leq a-n-2$ for $a\geq m_2+1$ we get
 \begin{align*}
\frac{|\Omega_{S_{j+1,i}^{a}}^{j+1}|}{|\Omega|}&\leq \sum_{l=1}^{i-1}C_{j,l}^{1}Cj^{2(m_1'+m_2'-n')}l^{m_1'}i^{2m_2'-a-n'}l^{-m_1'-2+\ve'}+C_{j,i}^{1}C(j^{2}i)^{m_1'+a-n'}i^{-m_1'-2+\ve'}\\
&\quad+\sum_{b=m_2+1}^{a-1}C_{j,i-1}^{2}C(i-1)^{b-a}(i+1)^{b-n'}(j+2-i)^{-2}\\
&\quad+C_{j,i-1}^{2}\left(\left(\frac{i+1}{i+2}\right)^{n'-a}+C(j(i-1))^{-2}\right)(i+1)^{a-n'}(j+2-i)^{-2}\\
&\quad+\sum_{b=a+1}^{n-m_1-1}C_{j,i-1}^{2}C(j^{2}(i-1))^{a-b}(i+1)^{b-n'}(j+2-i)^{-2}\\
&\leq \left(C_{j,i-1}^{2}+n4^{n}CC_{j,i-1}^{2}i^{-2}+4^{n}CM_{j}i^{-2+\ve'}\right)(i+2)^{a-n'}(j+2-i)^{-2}\\
&\leq C_{j+1,i}^{2}(i+2)^{a-n'}(j+2-i)^{-2}.
\end{align*}
 If  $\supp(\nu_{k})\cap S_{i,j+1}^{a}\neq \emptyset$, then
  \[\overline{\nu}_{k}=A_{k}\in \bigcup_{l=1}^{i}B_{l}^{j}\cup\bigcup_{a=m_2+1}^{n-m_1-1} S_{i-1,j}^{a},\]
  so, using again $2m_2-a-n\leq a-n-2$ for $a\geq m_2+1$ and that $a\leq n-m_1-1$ we obtain
\begin{align*}
\frac{|\Omega_{S_{i,j+1}^{a}}^{j+1}|}{|\Omega|}&\leq \sum_{l=1}^{i-1}C_{j,l}^{1}Cj^{m_1'+m_2'-n'}\left((i-1)^{2}(j+1)\right)^{m_1'+a-n'}l^{-2+\ve'}(j+2)^{m_2'-n'}\\
&\quad+C_{j,i}^{1}Cj^{m_2'-a}i^{-2+\ve'}(j+2)^{m_2'-n'}+\sum_{b=m_2+1}^{a-1}C_{j,i-1}^{2}Cj^{b-a}(j+2)^{b-n'}(j+2-i)^{-2}\\
&\quad+C_{j,i-1}^{2}\left(\left(\frac{j+2}{j+3}\right)^{n'-a}+C((i-1)j)^{-2}\right)(j+2)^{a-n'}(j+2-i)^{-2}\\
&\quad+\sum_{b=a+1}^{n-m_1-1}C_{j,i-1}^{2}C((i-1)^{2}j)^{a-b}(j+2)^{b-n'}(j+2-i)^{-2}\\
&\leq \left(C_{j,i-1}^{2}+n8^{n}CC_{j,i-1}^{2}i^{-2}+8^{n}CM_{j}i^{-2+\ve'}\right)(j+3)^{a-n'}(j+2-i)^{-2}\\
&\leq C_{j+1,i}^{2}(j+3)^{a-n'}(j+2-i)^{-2}.
\end{align*}
It only remains to estimate $|\Omega_{S_{j+1,j+1}^{a}}^{j+1}|$ for $a\in\{m_2+1,\ldots,n-m_1-1\}$.  If  $\supp(\nu_{k})\cap S_{j+1,j+1}^{a}\neq \emptyset$, then
  \[\overline{\nu}_{k}=A_{k}\in \bigcup_{l=1}^{j}\left(A_{j}^{l}\cup B_{l}^{j}\right)\cup\bigcup_{a=m_2+1}^{n-m_1-1} S_{j,j}^{a}.\]
Therefore, using $2m_2-a-n\leq a-n-2$ for $a\geq m_2+1$ we have
\begin{align*}
&\frac{|\Omega_{S_{j+1,j+1}^{a}}^{j+1}|}{|\Omega|}\leq \sum_{l=1}^{j}C_{j,l}^{1}C\left(j^{2(m_1'+m_2'-n')}l^{m_1'}(j+1)^{2m_2'-a-n'}l^{-m_1'-2+\ve'}\right.\\
&\quad+\left.j^{m_1'+m_2'-n'}\left(j^{2}(j+1)\right)^{m_1'+a-n'}l^{-2+\ve'}(j+2)^{m_2'-n'}\right)+\sum_{b=m_2+1}^{a-1}C_{j,j}^{2}Cj^{b-a}(j+2)^{b-n'}\\
&\quad+C_{j,j}^{2}\left(\left(\frac{j+2}{j+3}\right)^{n'-a}+C(j^{2})^{-2}\right)(j+2)^{a-n'}+\sum_{b=a+1}^{n-m_1-1}C_{j,j}^{2}C(j^{3})^{a-b}(j+2)^{b-n'}\\
&\leq \left(C_{j,j}^{2}+n8^{n}CC_{j,j}^{2}j^{-2}+8^{n}CM_{j}j^{-4}\right)(j+3)^{a-n'}\leq C_{j+1,j+1}^{2}(j+3)^{a-n'}.
\end{align*}
Hence, we have proved \ref{pro th ex border})-\ref{pro th ex Sij}). To see that there exists a convex function $u$ such that $f=\nabla u$ we observe that $Df\in\Gamma_{+}$ a.e. in $\Omega$ and we apply the reasoning of \cite[Th.1]{FaMoOl16} based on Poincar\'{e}'s lemma. Therefore, the proof of Theorem \ref{theorem existence} is finished.
\end{proof}
In fact, following the same reasoning as in \cite{LiMa16}, one can show that our function $u$ is strictly convex.

Note that the set $\bigcup_{i=1}^{j}\bigcup_{a=m_2+1}^{n-m_1-1} ( A_{j}^{i}\cup B_{i}^{j}\cup S_{j,i}^{a}\cup S_{i,j}^{a} )$ appearing in item \ref{pro th ex sup}) at the beginning of the proof of Theorem \ref{theorem existence} is the set $E^j_{m_1}$ mentioned in the Introduction, while $E^j_{m_2}$ is the inverse of that set.

\section*{Acknowledgements}

The author has been supported by the ERC Starting grant no.\ 307179 and this work is part of his doctoral thesis at Universidad Aut´onoma de Madrid under the supervision of Daniel Faraco and Carlos Mora-Corral. He would like to thank them for the suggestion of the problem and all the help provided during the preparation of this work.

{\small
\bibliography{Bibliography}
\bibliographystyle{siam}
}

\end{document}